\RequirePackage{ifpdf}
\ifpdf % We are running pdfTeX in pdf mode
\documentclass[pdftex]{sigma}
\else
\documentclass{sigma}
\fi

\numberwithin{equation}{section}
\newtheorem{Theorem}{Theorem}[section]
\newtheorem{Corollary}[Theorem]{Corollary}
\newtheorem{Lemma}[Theorem]{Lemma}
\newtheorem{Proposition}[Theorem]{Proposition}
{\theoremstyle{definition}
\newtheorem{Remark}[Theorem]{Remark}
}

\usepackage[thinlines]{easymat}
\usepackage{bm}

\newcommand\CC{{\mathbb C}}

\newcommand\RR{{\mathbb R}}

\newcommand\NN{{\mathbb N}}

\newcommand{\diag}{\operatorname{diag}}
\newcommand{\supp}{\operatorname{supp}}

\newcommand{\D}{\displaystyle}

\begin{document}

\allowdisplaybreaks

\renewcommand{\PaperNumber}{098}

\FirstPageHeading

\ShortArticleName{Properties of Matrix Orthogonal Polynomials via their Riemann--Hilbert Characterization}

\ArticleName{Properties of Matrix Orthogonal Polynomials\\ via their Riemann--Hilbert Characterization}

\Author{F.~Alberto GR\"UNBAUM~$^\dag$, Manuel D.~DE LA IGLESIA~$^\ddag$\\ and Andrei MART\'{I}NEZ-FINKELSHTEIN~$^\S$}

\AuthorNameForHeading{F.A.~Gr\"unbaum, M.D.~de la Iglesia and A.~Mart\'{i}nez-Finkelshtein}

\Address{$^\dag$~Department of Mathematics, University of California,  Berkeley, Berkeley, CA 94720  USA}
\EmailD{\href{mailto:grunbaum@math.berkeley.edu}{grunbaum@math.berkeley.edu}}
\URLaddressD{\url{http://math.berkeley.edu/~grunbaum/}}

\Address{$^\ddag$~Departamento de An\'{a}lisis Matem\'{a}tico,
   Universidad de Sevilla,\\
\hphantom{$^\ddag$}~Apdo (P. O. BOX) 1160, 41080 Sevilla, Spain}
\EmailD{\href{mailto:mdi29@us.es}{mdi29@us.es}}
\URLaddressD{\url{http://euler.us.es/~mdi29/}}

\Address{$\S$~Departamento de
Estad\'{i}stica y Matem\'{a}tica Aplicada, Universidad de Almer\'{i}a,\\
\hphantom{$^\S$}~04120
Almer\'{i}a, Spain}
\EmailD{\href{mailto:andrei@ual.es}{andrei@ual.es}}
\URLaddressD{\url{http://www.ual.es/~andrei/}}

\ArticleDates{Received June 09, 2011, in f\/inal form October 20, 2011;  Published online October 25, 2011}

\Abstract{We give a Riemann--Hilbert approach to the theory of matrix
orthogonal polynomials. We will focus on the algebraic aspects of the
problem, obtaining dif\/ference and dif\/ferential relations satisf\/ied
by the corresponding orthogonal polynomials. We will show that in the matrix case there is some extra freedom that allows us to obtain a family of ladder operators, some of them of 0-th order, something that is not possible in the scalar case. The combination of the ladder operators will lead to a family of second-order dif\/ferential equations satisf\/ied by the orthogonal polynomials, some of them of 0-th and f\/irst order, something also impossible in the scalar setting. This shows that the dif\/ferential properties in the matrix case are much more complicated than in the scalar situation. We will study several examples given in the last years as well as others not considered so far.}


\Classification{42C05; 35Q15}

\section{Introduction}

The theory of matrix orthogonal polynomials on the real line (MOPRL) has its foundations in the seminal papers of Krein~\cite{MR0034964, Krein:1971fk} (see also their account in the book of Berezans'ki{\u\i} \cite{MR0222718}). For further historical background and analytic results, the reader is referred to the survey~\cite{MR2379691} and to Chapter~4 of~\cite{Simon:2010fk}. In many aspects the MOPRL resemble their scalar counterparts, especially where the proofs are based on basic properties of Hilbert spaces. Nevertheless, the non-commutativity of matrix multiplication and the existence of non-zero singular matrices add features to the theory that make MOPRL an interesting object of study. Moreover, many problems for scalar polynomials are better understood or recast in terms of some matrix polynomials, see, for instance,~\cite{MR1327404}.

Recall that a matrix polynomial of degree $\leq n$ in $\CC^{N\times N}$  and a scalar variable $x$ can be def\/ined as an expression of the form
\begin{gather*}
\bm A_n x^n  + \cdots + \bm A_1 x + \bm A_0,
\end{gather*}
where $\bm A_j$'s are constant matrices in $\CC^{N\times N}$. In what follows we consider $N$ f\/ixed, and denote by $\mathbb P_n$ the family of matrix polynomials of degree $\leq n$ in $\CC^{N\times N}$, as well as by $\mathbb P := \bigcup_{n\geq 0} \mathbb P_n$. We use preferably boldface letters to denote matrices, and standard font for scalars.
We  also use $\bm I_N$ for the $N\times N$ identity matrix, omitting the explicit reference to its dimension when it cannot lead anyone into confusion.

In this paper we consider some aspects of the MOPRL theory, assuming that the orthogonality is given by an absolutely continuous measure on the line. More precisely, our starting point is a weight $\bm W=(W_{ij}): (a,b) \to GL(N, \mathbb R)$, def\/ined and positive def\/inite on a f\/inite or inf\/inite interval $(a,b)\subset \mathbb R$. We will assume that  all $W_{ij}$ and $W'_{ij}$  have f\/inite moments:
\begin{gather*}
\int_a ^b |x|^n \bm W(x) dx <\infty, \qquad \int_a ^b |x|^n \bm W'(x) dx <\infty,  \qquad n \in \mathbb N_0 := \mathbb N \cup \{0\},
\end{gather*}
where the integration of a matrix function is applied entry-wise.

For any two $\bm P, \bm Q\in \mathbb P$, the weight $\bm W$ induces two matrix-valued ``inner products'',
\begin{gather*}\label{Innerp}
    (\bm P,\bm Q)_{\bm W}=\int_a ^b \bm P(x)\bm W(x)\bm Q^*(x)\, dx,
\end{gather*}
and
\begin{gather*}
 \langle \bm P,\bm Q\rangle_{\bm W}=\int_a ^b  \bm Q^*(x) \bm W(x) \bm P(x) \, dx =(\bm Q^*, \bm P^*)_{\bm W},
\end{gather*}
where the asterisk denotes the conjugate transpose  (or Hermitian conjugate) of a matrix. Due to this connection between both inner products, we restrict our attention to $(\cdot, \cdot)_{\bm W}$. We def\/ine also the norm
\begin{gather*}
\| \bm P\|_W := \left( \text{Tr}\,  \langle \bm P,\bm P\rangle_{\bm W} \right)^{1/2},
\end{gather*}
and assume that $\bm W$ is non-trivial, in the sense that $\|\bm P\|_{\bm W}>0$ for every non-zero matrix polynomial $\bm P$. In this case (see \cite[Lemma~2.3]{MR2379691} or \cite[Proposition~4.2.3]{Simon:2010fk}), $(\bm P,\bm P)_{\bm W}$ is non-singular for every non-zero polynomial $\bm P$, and we can easily implement a matrix analogue of the Gram--Schmidt orthogonalization procedure, which yields a unique sequence $(\widehat{\bm P}_n)_n$ of monic orthogonal polynomials such that $\widehat{\bm P}_0=\bm I$,
\begin{gather}
\label{orthogonality}
\widehat{\bm P}_n(x)=x^n \bm I +\sum_{j=0}^{n-1} \bm a_{n,j}x^j, \qquad (\widehat{\bm P}_n, \bm Q)_{\bm W}=\bm{0} \quad \text{for every} \quad \bm Q\in \mathbb P_{n-1}, \quad n \in \NN,
\end{gather}
as well as matrix polynomials $(\bm P_n)_n$, ``orthonormal''  with respect to $\bm W$, such that
\begin{gather} \label{def:leading}
\bm P_n(x)=\bm \kappa_n \widehat{\bm P}_n, \qquad  (\bm P_n,\bm P_m)_{\bm W}= \delta_{n,m}\bm I.
\end{gather}
Obviously, $\bm   P_n$'s are determined up to a unitary left factor, so we can speak about an equivalence class of orthonormal MOPRL, corresponding to the weight $\bm W$. Additionally, if $\bm B$ is a constant non-singular matrix, then
\begin{gather*}
\widetilde{\bm W}(x)=\bm B \bm W(x) \bm B^*
\end{gather*}
is also a weight, and $(\widetilde{\bm P}_n)_n=(\bm P_n\bm B^{-1})_n$ is the corresponding sequence of orthonormal polynomials. Hence, more than a single weight we consider an equivalence class of weights given by its representative $\bm W$. In particular, without loss of generality we can assume that $\bm W(x_{0})=\bm I$ at a point $x_{0}\in (a,b)$. If this class contains a diagonal matrix-valued function, we say that $\bm W$ reduces to \emph{scalar weights}. A characterization of this fact, as shown in~\cite{MR2039133}, is the commutativity condition
\begin{gather*}
\bm W(x) \bm W(y)=\bm W(y) \bm W(x), \qquad x, y \in [a,b].
\end{gather*}
This situation is considered trivial, and usually omitted from consideration.

Work in the last few years has revealed a number of explicit families of MOPRL; in many cases they are joint eigenfunctions of some f\/ixed dif\/ferential operator with matrix coef\/f\/icients independent of the degree $n$ of the MOPRL. This study was initiated in \cite{MR1466158}, but nontrivial examples had to wait until \cite{MR2039133, MR1988656, MR1883412, MR2075953}. A solution to the classif\/ication problem for those $\bm W$ whose MOPRL are common eigenfunctions of some f\/ixed dif\/ferential operator remains elusive. There are by now two methods that have yielded nontrivial examples: a) the connection with matrix-valued spherical functions for symmetric spaces~\cite{MR1883412, MR2075953}, and b) a combination of classical methods and some Lie algebra tools \cite{MR2039133, MR2152234}. A necessary condition in terms of moments is given in~\cite{MR2039133}, and a necessary and suf\/f\/icient condition in terms of the ``ad-conditions'' is given in \cite{GT}. In this approach one gets an explicit formula for one dif\/ferential operator. Concerning applications, the recurrent relations have been used in the study of certain quasi-birth-and-death processes~\cite{MR2407600, MR2568634, MR2723279, MR2421469}, in which case the dif\/ferential operator plays no role. Recently, new applications of these processes have been related to urn or Young diagrams models \cite{GPT}.

A big boost in the research of the scalar OPRL has been their Riemann--Hilbert (RH) cha\-racterization, introduced in the seminal paper of Fokas, Its and Kitaev~\cite{Fokas92}, and complemented with a non-linear steepest descent analysis in a series of works of Deift, Zhou and collaborators~\cite{MR2000g:47048, MR2307753, MR94d:35143, MR1989226}. This combination has allowed to establish extremely strong asymptotic results with applications in  random matrix theory, in particular for unitary invariant ensembles~\cite{MR2000g:47048, MR2514781}, determinantal point processes~\cite{MR2327475, MR2450071, MR2470930, MR2647568}, orthogonal Laurent polynomials~\cite{MR2253536, MR2336421}, Painlev\'e trascendents~\cite{MR2488326, MR2519676}, the Toda lattices~\cite{MR1369835}, to mention a few.

The RH characterization, even without the steepest descent analysis, allows one to prove other algebraic and analytic properties of orthogonal polynomials, as it was illustrated in~\cite{MR2307753} for OP on the unit circle, and in~\cite[Chapter 22]{Ismail05}, for classical families of OPRL.  One of the goals of this paper is to extend these considerations to MOPRL, something that, surprisingly enough, has not been explored in depth so far. We will show that the RH formulation, which has a~very natural block-wise generalization to the MOPRL case, allows one  to reveal some subtleties hidden in the matrix case. For instance, we will derive identities and dif\/ference-dif\/ferential equations that were unknown even for some of the explicit families mentioned above. The RH technique can be also a natural method for obtaining the dif\/ferential relations for MOPRL with respect to general weight matrices.

There have been several papers considering the block-wise RH problem of the f\/lavor similar to the one discussed in this paper, but for singular (namely, rank 1) weight matrices. For these weights the MOPRL are connected with a class of multiple orthogonal polynomials that f\/ind applications in the analysis of determinantal point processes and non-intersecting stochastic paths (see e.g.\ the works of Kuijlaars and collaborators \cite{MR2327475, MR2610344, MR2539187, MR2647568}) or of the multicomponent 2D Toda lattice hierarchy \cite{2010InvPr..26e5009A}.

This paper deals with strictly algebraic consequences of the RH formulation of matrix polynomials orthogonal with respect to a weight matrix supported on an interval $[a, b]\subset\mathbb{R}$; the asymptotic analysis via the Deift--Zhou non-linear steepest descent method will not be discussed here.

In Section~\ref{sec:1} we will discuss two dual RHPs, which are uniquely solved in terms of the MOPRL. Standard arguments related to RHP allow us to obtain the three-term recurrence relation and the dif\/ferential identities (ladder operators) for MOPRL. Both yield also the so-called Lax pair for the polynomials, which is an over-determined system~\eqref{LaxP}, whose compatibility conditions render further properties of MOPRL. All the results in this section have been previously obtained using dif\/ferent approaches, but the RH method will provide new and straightforward proofs.

In Section~\ref{secRH} we will consider a special transformation of the original RHP, based on facto\-ri\-zations of the weight matrix of the form $\bm W=\bm T\bm T^*$. The main goal is to obtain a RHP with constant jumps (independent of $z$), in order to simplify the dif\/ferential relation obtained earlier. For simplicity, we will focus on the case where the support of the weight matrix is $\mathbb{R}$. The extension of the method to weights supported on semi-inf\/inite or f\/inite intervals is simple, see a~short discussion at the end of this paper.

We exploit also the non-uniqueness of the weight factorization, which in the matrix case gives us some extra freedom and yields non-trivial relations. In particular, a family of ladder operators, some of them of 0-th order, is obtained, a phenomenon that is not possible in the scalar situation. Furthermore, by combining appropriately the family of ladder operators, we will get the 0-th, f\/irst and second order dif\/ferential equations satisf\/ied by the MOPRL. Some examples of lower order dif\/ferential operators appear in \cite{MR2217096, MR2357983, MR2083080}.

The considerations so far have been completely general, applicable to any family of MOPRL. In Section~\ref{sec:examples} we narrow the analysis of Sections~\ref{sec:1} and~\ref{secRH} to some relevant examples of MOPRL  supported on $\mathbb{R}$, obtaining new results even in the cases  studied previously in the literature.

In the f\/inal stage of preparation of this manuscript we learned about the preprint \cite{CCM} which also uses the Riemann--Hilbert approach for the analysis of the matrix orthogonal polynomial. Although there is some overlapping between the results contained in our Sections~\ref{sec:1},~\ref{secRH} and in~\cite[Section~2]{CCM}, the focus of both contributions is dif\/ferent, and in this sense, complementary.

\section[The Riemann-Hilbert problem for MOPRL]{The Riemann--Hilbert problem for MOPRL} \label{sec:1}

\subsection{Formulation and basic properties} \label{sec:RHformulation}

In this section we discuss the Riemann--Hilbert problem (RHP) related to MOPRL with respect to a $N \times N$ weight matrix $\bm W$ supported on an interval $[a\,, b\, ]\subset\mathbb{R}$; this interval can be either bounded or unbounded. We assume $\bm W$ continuous and non-vanishing on $(a,b)$, and that at any f\/inite endpoint of the support the weight $\bm W$ has at worse a power-type singularity, that is, $\bm W$ is of the form
\begin{gather}\label{BouW}
   \bm W(z)=|z-c|^{\gamma_c}\widetilde{\bm W}(z), \qquad  \gamma_c>-1,
\end{gather}
where $c\in \{a, b\}$, $c\neq \pm \infty$, and $\widetilde{\bm W}$ is a bounded, continuous and non-vanishing at $z=c$. This class of weights comprises all the examples considered so far in the literature.

As a convention, in what follows we write the $2N \times 2N$ matrices partitioned into $N\times N$ blocks, as in \eqref{jump_cond} below.  Recall that $\bm A^*$ stands for the Hermitian conjugate of the matrix $\bm A$, as well as $\bm A^{-*}=(\bm A^*)^{-1}$.  We adopt the convention that for any  matrix-valued function $\bm P(z)$ of a complex variable $z$,
$\bm P^*$ denotes the matrix-valued function obtained as
\begin{gather*}
    \bm P^*(z) := \left( \bm P(\bar z) \right)^*.
\end{gather*}

The RHP for MOPRL with respect to a weight matrix $\bm W$ consists in f\/inding a matrix function $\bm Y^{n}:\mathbb{C}\rightarrow\mathbb{C}^{2N\times 2N}$ such that
\begin{enumerate}\itemsep=0pt
        \item[(Y1)] $\bm Y^{n}$ is analytic in $\mathbb{C}\setminus[a, b]$.
        \item[(Y2)] $\bm Y^{n}$ has on $(a,b)$ continuous boundary values $\bm Y^{n}_+$ (resp., $\bm Y^{n}_-$) from the upper (resp., lower) half plane, such that
      \begin{gather}\label{jump_cond}
        \bm Y^{n}_{+}(x)=\bm Y^{n}_-(x) \left(\begin{MAT}{c.c}
  \bm I_N & \bm W(x) \\.
 \bm  0 & \bm I_N \\
\end{MAT} \right), \qquad x\in(a, b).
      \end{gather}
        \item[(Y3)] As $z\rightarrow\infty$, for every $m\in \mathbb N$ we have
\begin{gather}\label{Yninf}
\bm Y^n(z)=\left(\bm I_{2N}+\sum_{i=1}^m \dfrac{\bm Y_i^n}{z^i}+\mathcal{O} \big(1/z^{m+1}\big)\right)
\left(\begin{MAT}{c.c}
  z^n\bm I_N & \bm 0 \\.
        \bm  0 & z^{-n}\bm I_N \\
\end{MAT} \right)
\end{gather}
(where the asymptotic term $\mathcal{O} (1/z^{m+1})$ depends on $n$).
        \item[(Y4)] As $z\rightarrow c$, $c\in \{a, b\}$, $c\neq \pm \infty$, we have
       \begin{gather*}
       \bm Y^{n}(z)=        \left(\begin{MAT}{c.c}
 \mathcal{O}(1) & \mathcal{O}(h(z)) \\.
  \mathcal{O}(1) & \mathcal{O}(h(z))\\
\end{MAT} \right),
       \end{gather*}
       where
      \begin{gather}
      \label{def_h}
      h(z)=\begin{cases}
      |z-c|^{\gamma_c}, & \text{if } -1<\gamma_c<0, \\
      \log |z-c| , & \text{if }  \gamma_c=0, \\
      1, & \text{if }  \gamma_c>0.
      \end{cases}
      \end{gather}
\end{enumerate}
\begin{Remark}
Usually the asymptotic condition at inf\/inity \eqref{Yninf} is stated for $m=1$; however, it can be proved that \eqref{Yninf} for $m=1$ implies that this condition holds for every $m\in \NN$.
\end{Remark}

Along with the RH problem (Y1)--(Y4) we can consider the following dual problem: f\/inding a matrix function $\bm y^{n}:\mathbb{C}\rightarrow\mathbb{C}^{2N\times 2N}$ such that
\begin{enumerate}\itemsep=0pt
        \item[(y1)] $ \bm y^{n} $ is analytic in $\mathbb{C}\setminus[a, b]$.
        \item[(y2)] $ \bm y^{n}_{+}(x)= \left(\begin{MAT}{c.c}
  \bm I_N & -\bm W(x) \\.
 \bm  0 & \bm I_N \\
\end{MAT} \right)\bm y^{n}_-(x)$ when $x\in(a, b)$.
        \item[(y3)] As $z\rightarrow\infty$, for every $m\in \mathbb N$ we have
\begin{gather}\label{Yninfi}
\bm y^{n}(z)= \left(\begin{MAT}{c.c}
z^{-n}\bm I_N & \bm 0 \\.
 \bm 0 & z^n\bm I_N \\
\end{MAT} \right) \left(\bm I_{2N}+ \sum_{i=1}^m \dfrac{\widetilde{\bm Y}_i^n}{z^i}+\mathcal{O} \big(1/z^{m+1}\big)\right).
\end{gather}
        \item[(y4)] As $z\rightarrow c$, $c\in \{a, b\}$, $c\neq \pm \infty$,  we have
\begin{gather*}\bm y^{n}(z)=  \left(\begin{MAT}{c.c}
 \mathcal{O}(h(z)) & \mathcal{O}(h(z)) \\.
  \mathcal{O}(1) & \mathcal{O}(1)\\
\end{MAT} \right),
 \end{gather*}
 with $h$ def\/ined in \eqref{def_h}.
\end{enumerate}
It will turn out that (y1)--(y4) is related to the inverse of the solution of (Y1)--(Y4).

For any integrable $N\times N$ matrix-valued function $\bm F$ on $[a,b]$,
\begin{gather*}
\mathcal C(\bm F)(z) := \frac{1}{2\pi i}  \int_a^b  \dfrac{\bm F(t)}{t-z} dt
\end{gather*}
def\/ines the Cauchy or Stieltjes transform of $\bm F$, which is a matrix-valued and analytic function in $\CC\setminus [a,b]$. Let us introduce the matrix polynomials of the second kind (of degree $n-1$),
def\/ined by
\begin{gather}\label{QQn}
   \bm  Q_n(x)=\int_a^b  \dfrac{\bm P_n(t)-\bm P_n(x)}{t-x}\bm W(t)dt,\qquad
    n\geq0.
\end{gather}
We have
\begin{gather*}
2\pi i\mathcal C(\bm P_n\bm W)(z)=\bm Q_n(z)+2\pi i \bm P_n(z)\mathcal C(\bm W)(z),
\end{gather*}
where $2\pi i \mathcal C(\bm W)$ is the $m$-function of the weight matrix $\bm W$.

Finally, if  $\bm \kappa_{n}$ is the leading coef\/f\/icient of any corresponding normalized polynomial $\bm P_n$, we denote
\begin{gather}
\label{def:gamma}
\bm \gamma_n=\bm \kappa_{n}^* \bm \kappa_{n}.
\end{gather}
Observe that a priori $\bm \gamma_n$ depends on the selection of $\bm \kappa_{n}$.

The following is a complete analogue of the well-known theorem from \cite{Fokas92}, although its proof requires some additional considerations:
\begin{Theorem}\label{thm}
The unique solution of the RHP {\rm (Y1)--(Y4)} is
\begin{gather}\label{RHS}
\bm Y^{n}(z)=\bm R_n(z)  \bm Y^{0}(z),\qquad n\geq 0,
\end{gather}
where
\begin{gather*}
\bm Y^{0}(z)= \left(\begin{MAT}{c.c}
  \bm I_N & \mathcal C(\bm W)(z) \\.
   \bm 0 & \bm I_N\\
\end{MAT} \right),
\end{gather*}
and the transfer matrix $\bm R_n$ is a matrix polynomial given by $\bm R_{0}(z)=\bm I$,
\begin{gather}\label{TransferM}
 \bm R_{n}(z)= \left(\begin{MAT}{c.c}
 \bm \kappa_n^{-1}\bm P_n(z) & ( 2\pi i\bm \kappa_n)^{-1} \bm Q_n(z) \\.
 -2\pi i\bm \kappa_{n-1}^* \bm P_{n-1}(z) & -\bm \kappa_{n-1}^* \bm Q_{n-1}(z) \\
\end{MAT} \right), \qquad n\in \NN.
\end{gather}
Analogously, the unique solution of the RHP {\rm (y1)--(y4)} is
\begin{gather}\label{rhsy}
\bm y^{n}(z)=\bm y^{0}(z)  \bm r_n(z),\qquad n\geq 0,
\end{gather}
where
\begin{gather*}
\bm y^{0}(z)=  \left(\begin{MAT}{c.c}
  \bm I_N & -\mathcal C(\bm W)(z) \\.
   \bm 0 & \bm I_N\\
\end{MAT} \right),
\end{gather*}
and the transfer matrix $\bm r_n$ is a matrix polynomial given by $\bm r_{0}(z)=\bm I$,
\begin{gather}\label{Transferr}
 \bm r_{n}(z)=\left(\begin{MAT}{c.c}
   -\bm Q_{n-1}^*(z)\bm \kappa_{n-1} & -( 2\pi i)^{-1} \bm Q_n^*(z)\bm \kappa_n^{-*} \\.
  2\pi i \bm P_{n-1}^*(z) \bm \kappa_{n-1} & \bm P_n^*(z) \bm \kappa_n^{-*} \\
\end{MAT} \right), \qquad n\in \NN.
\end{gather}

Moreover, for all $n\geq 0$,
\begin{gather}
%\begin{split}
\label{connection}
  \bm y^n(z)  =\left(\begin{MAT}{c.c}
 \bm 0 & -\bm I_N \\.
 \bm I_N &  \bm 0 \\
\end{MAT} \right) (\bm Y^n(z))^{T}\left(\begin{MAT}{c.c}
 \bm 0 &  \bm I_N \\.
 -\bm I_N &  \bm 0 \\
\end{MAT} \right) =(\bm Y^n(z))^{-1}, \\   \det \bm Y^{n}(z)  =1, \qquad \text{for all}  \ \ z\in\mathbb{C}.
\label{det1}
\end{gather}
Here $\bm A^T$ denotes the transpose of the matrix $\bm A$.
\end{Theorem}
\begin{Remark}
In the scalar case, when $\bm Y^n$ is a $2\times 2$ matrix, there is no need to consider simultaneously both RHPs (Y1)--(Y4) and (y1)--(y4) since the existence of $(\bm Y^n)^{-1}$ and~\eqref{Transferr},~\eqref{connection} follow directly from \eqref{det1}. In the general $(2N)\times (2N)$ case more care should be put in the analysis of the local behavior at the endpoints $a$, $b$,  and \eqref{connection} is no longer straightforward.
\end{Remark}

\begin{Remark}
The solution $\bm Y^n$ of (Y1)--(Y4) satisf\/ies the following symmetry relation:
\begin{gather}\label{symmetryForYnew}
\bm Y^n(z)  =\left(\begin{MAT}{c.c}
  \bm I_N & \bm 0 \\.
\bm 0 & - \bm I_N \\
\end{MAT} \right) \overline{\bm Y^n(\bar{z})}\left(\begin{MAT}{c.c}
   \bm I_N & \bm 0 \\.
\bm 0 & -\bm I_N  \\
\end{MAT} \right),
\end{gather}
which yields some obvious consequences for the block entries of $\bm Y^n$. This relation is established using the invariance of (Y1)--(Y4) by such a conjugation.
\end{Remark}
\begin{Remark}
Alternatively, we can write the solution \eqref{RHS} as
\begin{gather}\label{RHS1}
\bm Y^{n}(z)=\left(\begin{MAT}{c.c}
 \widehat{\bm P}_n(z) & \mathcal C(\widehat{\bm P}_n\bm W)(z) \\.
 -2\pi i \bm \gamma_{n-1}\widehat{\bm P}_{n-1}(z) & -2\pi i\bm \gamma_{n-1}\mathcal C(\widehat{\bm P}_{n-1}\bm W)(z)  \\
\end{MAT} \right),\qquad  n\in \NN,
\end{gather}
where $\widehat{\bm P}_n$ denote the monic MOPRL of degree $n$ for the weight matrix $\bm W$. In the same vein,
\begin{gather}\label{RHSi}
(\bm Y^{n})^{-1}(z)=\left(\begin{MAT}{c.c}
  -2\pi i\mathcal C(\bm W\widehat{\bm P}^*_{n-1})(z)\bm \gamma_{n-1}& -\mathcal C(\bm W\widehat{\bm P}^*_n)(z) \\.
  2\pi i\widehat{\bm P}^*_{n-1}(z)\bm \gamma_{n-1} & \widehat{\bm P}^*_n(z) \\
\end{MAT} \right),\qquad n\in \NN.
\end{gather}
\end{Remark}

\begin{Remark}
Although $\bm \kappa_{n}$ is def\/ined up to a left unitary factor, the matrix coef\/f\/icient~$\bm \gamma_{n}$ in~(\ref{def:gamma}) is unique. This is a consequence of the uniqueness of the solutions of the RHP above.
\end{Remark}

\begin{Remark}
Formulas \eqref{RHS} and \eqref{rhsy} show that generically the behavior of $\bm Y^{n}$ and $(\bm Y^{n})^{-1}$ at the endpoints (and in general, any singular point) of the support of $\bm W$ is given by the local behavior of the Cauchy transform of the orthogonality weight, also known as its $m$-function, see \cite[\S~1.2]{Simon:2005}.
\end{Remark}

\begin{proof}[Proof of Theorem \ref{thm}]
The fact that \eqref{RHS} or \eqref{RHS1} is a solution of the RHP (Y1)--(Y4)  is established following the proof for the scalar case, see e.g.~\cite{MR2000g:47048}, and taking into account that  the orthogonality of $\widehat{\bm P}_n$ in \eqref{orthogonality} is equivalent to the homogeneous system
\begin{gather*}
\int_a^b x^j \widehat{\bm P}_n(x) \bm W(x)  \, dx =\bm 0, \qquad j=0, 1, \dots, n-1.
\end{gather*}
The same applies to \eqref{rhsy} or \eqref{RHSi} and the RHP (y1)--(y4).

Consider the function $\bm Z^n(z)= \bm Y^n(z) \bm y^n(z)$; from (Y2) and (y2) it follows that it has no jump across $(a,b)$, and by (Y4), (y4),
\begin{gather*}
\bm Z^n(z)=  \mathcal{O}(h(z)), \qquad z\to c\in\{a, b\}.
\end{gather*}
Hence, $\bm Z^n$ has only removable singularities at the f\/inite endpoints of the support of the weight, and thus is an entire function. It remains to observe that $\bm Z^n(\infty)=\bm I_{2N}$ to conclude that $\bm Z^n(z)=\bm I_{2N}$ for all $z\in \CC$, which proves that $ \bm y^n(z)=(\bm Y^n(z))^{-1}$, as well as the uniqueness of both solutions. The f\/irst identity in \eqref{connection} can be established by direct calculation or by observing that this transformation carries the RHP (Y1)--(Y4) to the RHP (y1)--(y4).

Finally, the scalar function $ \det\bm Y^n (z)$ is analytic across $[a,b]$, and by  \eqref{RHS}, \eqref{TransferM},
\begin{gather*}
\det \bm Y^{n}(z)=\det \bm R_{n}(z)
\end{gather*}
can have only removable singularities at the (f\/inite) endpoints of the support of $\bm W$. Hence, $ \det\bm Y^n (z)$ is an entire function; since by \eqref{Yninf}, $ \det\bm Y^n (\infty)=1$, we conclude that it is iden\-ti\-cally~1.
\end{proof}

\begin{Proposition}
Let $\bm P_n = \bm \kappa_n \widehat{\bm P}_n$, $n\geq 0$,  be a sequence of orthonormal MOPRL, and let  $(\bm Q_n)_n$ be the corresponding matrix polynomials of the second kind \eqref{QQn}. For these polynomials, define
\begin{gather}\label{An}
\bm A_n=\bm \kappa_{n-1}\bm \kappa_n^{-1}.
\end{gather}
Then
\begin{gather}\label{LOF}
 \bm  Q_n(z) \bm P_{n-1}^*(z)- \bm P_{n}(z) \bm Q_{n-1}^*(z)=\bm A_n^{-1},
\\
\label{HP}
  \bm   Q_n(z)\bm P_{n}^*(z)=\bm P_{n}(z)\bm Q_{n}^*(z).
\end{gather}
Moreover, $\bm P_{n-1}^*(x) \bm A_n \bm P_n(x)$ and $\bm Q_{n-1}^*(x) \bm A_n \bm Q_n(x)$ are Hermitian for all $x\in\mathbb{R}$ and $n\geq0$.
\end{Proposition}
\begin{Remark}
Identities \eqref{LOF} and \eqref{HP} are also known as the \emph{Liouville--Ostrogradski
formula} and the \emph{Hermitian property}, respectively. Both were originally derived for MOPRL in~\cite{MR1426899}, using a dif\/ferent approach.
\end{Remark}

\begin{proof}
Using the explicit expressions for $\bm Y^{n} $ and $(\bm Y^{n})^{-1}$ (written in the form \eqref{RHS} and \eqref{rhsy}), from block entries $(1,1)$ and $(2,2)$ of the identity $ \bm Y^{n} (\bm Y^{n})^{-1}=\bm I$ we obtain the Liouville--Ostrogradski formula \eqref{LOF} while from block entries $(1,2)$ and $(2,1)$ we get the so-called Hermitian property \eqref{HP}.

Block entries $(1,2)$ and $(2,1)$ of $(\bm Y^{n})^{-1} \bm Y^{n}=\bm I$ yield the
commutativity relations
\begin{gather*}\label{Crel}
 \bm P_{n-1}^*(z) \bm A_n \bm P_n(z)= \bm P_{n}^*(z) \bm A_n^* \bm P_{n-1}(z),\qquad
 \bm Q_{n-1}^*(z) \bm A_n \bm Q_n(z)= \bm Q_{n}^*(z) \bm A_n^* \bm Q_{n-1}(z),
\end{gather*}
which show that $\bm P_{n-1}^*(x) \bm A_n \bm P_n(x)$ and $\bm Q_{n-1}^*(x) \bm A_n \bm Q_n(x)$ are Hermitian for all $x\in\mathbb{R}$ and $n\geq0$.
\end{proof}

For the formulation of the following results we need to introduce a new set of parameters,
\begin{gather}\label{bnk}
   \bm b_{n,k}=\int_{a}^bx^k \widehat{\bm P}_n(x)\bm W(x)\, dx,\qquad k=n,n+1,\ldots,
\end{gather}
where $\widehat{\bm P}_n$ are the monic MOPRL. In the spirit of  \cite{MR2168386}, we can express them in terms of the coef\/f\/icients $\bm a_{n,j}$ of the polynomials $\widehat{\bm P}_n$ (see \eqref{orthogonality}) in a form suitable for numerical implementation:
\begin{Lemma}\label{lemmabn}
Let $\bm \Omega$ be the block lower triangular matrix built up from the coefficients of the MOPRL $(\widehat{\bm P}_n)_n$,
\begin{gather*}
\bm \Omega=\begin{pmatrix}
                \bm I  &  &   &  \\
                \bm a_{1,0} & \bm I  &   &  \\
                \vdots & \vdots & \ddots & \\
                 \bm a_{n,0}& \bm a_{n,1} &  \cdots & \bm I  \\
              \end{pmatrix},
\qquad \text{so that} \quad
\begin{pmatrix}
  \widehat{\bm P}_0(x) \\
  \widehat{\bm P}_{1}(x) \\
  \vdots \\
  \widehat{\bm P}_n(x) \\
\end{pmatrix}=\bm \Omega
\begin{pmatrix}
  \bm I \\
  x\bm I \\
  \vdots \\
  x^n\bm I \\
\end{pmatrix}, \qquad  \bm I= \bm I_N.
\end{gather*}
Then $\bm b_{n,k}$ defined in \eqref{bnk} can be obtained from the last $N\times N$ block row of  $\bm \Omega^{-1}$ as follows:
\begin{gather}\label{bnkx}
    \bm b_{n-k,n}^*\bm \gamma_{n-k}=\big(\bm \Omega^{-1}\big)_{n+1,n-k+1},\qquad k=0,\ldots,n.
\end{gather}
\end{Lemma}

\begin{proof}
Observe that
\begin{gather*}
\left(   \bm b_{0,n}^* ,  \dots , \bm b_{n-1,n}^*\bm ,    \bm b_{n,n}^*\bm  \right)   =\int_{a}^b x^n \bm W \big(
         \widehat{\bm P}^*_{0},  \widehat{\bm P}^*_{1},   \dots  , \widehat{\bm P}^*_{n}\big)dx \\
\hphantom{\left(   \bm b_{0,n}^* ,  \dots , \bm b_{n-1,n}^*\bm ,    \bm b_{n,n}^*\bm  \right) }{}
= \int_{a}^b x^n \bm W    \left(
           \bm I , x\bm I ,   \dots ,  x^{n}\bm I
         \right) \bm \Omega^*\, dx
 = \left( \bm \mu_{n},  \bm \mu_{n+1},  \dots  ,  \bm \mu_{2n} \right) \bm \Omega^*,
\end{gather*}
where $\bm \mu_{n}=\int_{a}^bx^n\bm W(x)dx$ are the moments of the weight matrix $\bm W$. Analogously, by \eqref{def:leading},
\begin{gather*}
\left( \diag \left( \bm \gamma_{k} \right)_{k=0}^n\right)^{-1}   = \diag \left( \bm \gamma_{k}^{-1} \right)_{k=0}^n   = \int_{a}^b\begin{pmatrix}
  \widehat{\bm P}_0 \\
  \widehat{\bm P}_{1} \\
  \vdots \\
  \widehat{\bm P}_n
\end{pmatrix}\bm W\big(
          \widehat{\bm P}^*_{0},  \widehat{\bm P}^*_{1},   \dots , \widehat{\bm P}^*_{n}
         \big)\, dx \\
\qquad{}
 = \bm \Omega\int_{a}^b\begin{pmatrix}
  \bm I \\
  x\bm I \\
  \vdots \\
  x^n\bm I \\
\end{pmatrix}\bm W\big(
           \bm I,x\bm I,  \dots,  x^{n}\bm I
         \big)\bm \Omega^*
          =\bm \Omega \begin{pmatrix}
                        \bm \mu_{0} & \bm \mu_{1} & \cdots & \bm \mu_{n} \\
                        \bm \mu_{1} & \bm \mu_{2}& \cdots & \bm \mu_{n+1} \\
                        \vdots & \vdots & \ddots & \vdots \\
                        \bm \mu_{n}& \bm \mu_{n+1}&\cdots& \bm \mu_{2n}
                      \end{pmatrix} \bm \Omega^*.
\end{gather*}
Hence,
\begin{gather*}
 \big(
         \bm b_{0,n}^*\bm \gamma_0,   \dots, \bm b_{n-1,n}^*\bm \gamma_{n-1},    \bm b_{n,n}^*\bm \gamma_n
         \big) \\
 \qquad{} =\big(
           \bm \mu_{n},  \bm \mu_{n+1},  \dots ,  \bm \mu_{2n}
         \big)\begin{pmatrix}
                        \bm \mu_{0} & \bm \mu_{1} & \cdots & \bm \mu_{n} \\
                        \bm \mu_{1} & \bm \mu_{2}& \cdots & \bm \mu_{n+1} \\
                        \vdots & \vdots & \ddots & \vdots \\
                        \bm \mu_{n}& \bm \mu_{n+1}&\cdots& \bm \mu_{2n}
                      \end{pmatrix}^{-1}\bm \Omega^{-1}
                        =\big(
          \bm 0,   \bm 0,   \dots, \bm I
         \big) \bm \Omega^{-1},
\end{gather*}
which yields \eqref{bnkx}.
\end{proof}

For what follows it is useful to single out the explicit expressions for  $\bm b_{n-3,n}, \dots,  \bm b_{n,n}$, that are obtained from Lemma \ref{lemmabn} by direct computations:
\begin{Corollary}\label{corobn}
Coefficients $\bm b_{n-k,n}$, for $k=0,1,2,3$, are given by
\begin{gather*}\bm \gamma_{n}\bm b_{n,n}=\bm I,\\
 \bm \gamma_{n-1}\bm b_{n-1,n} =-\bm a_{n,n-1}^*,\\
 \bm \gamma_{n-2}\bm b_{n-2,n}=\bm a_{n-1,n-2}^*\bm a_{n,n-1}^*-\bm a_{n,n-2}^*,\\
\bm \gamma_{n-3}\bm b_{n-3,n}=-\bm a_{n-2,n-3}^*\bm a_{n-1,n-2}^*\bm a_{n,n-1}^*+\bm a_{n-1,n-3}^*\bm a_{n,n-1}^*+\bm a_{n-2,n-3}^*\bm a_{n,n-2}^*-\bm a_{n,n-3}^*.
\end{gather*}
\end{Corollary}

Now we return to Theorem \ref{thm}; as its immediate consequence we can relate  the coef\/f\/icients in the asymptotic expansion of $\bm Y ^n$ and of $(\bm Y ^n)^{-1}$, with the coef\/f\/icients $\bm a_{n,j}$ and $\bm b_{n,j}$:

\begin{Corollary}\label{coro25}
The coefficients $\bm Y_i^n$ in  \eqref{Yninf} are given by
\begin{gather}\label{Yni}
   \bm Y^n_i=
   \left(\begin{MAT}{c.c}
  \bm a_{n,n-i} & -\D\frac{1}{2\pi i}\bm b_{n,n+i-1} \\.
 -2\pi i\bm \gamma_{n-1}\bm a_{n-1,n-i} & \bm \gamma_{n-1}\bm b_{n-1,n+i-1} \\
\end{MAT} \right),\qquad i\geq0.
\end{gather}

Analogously, the coefficients $\widetilde{\bm Y}_i^n$ in \eqref{Yninfi} are given by
\begin{gather}\label{Ynig}
    \widetilde{\bm Y}^n_{i}=\left(\begin{MAT}{c.c}
            \bm b_{n-1,n+i-1}^T \bm \gamma_{n-1} & \D\frac{1}{2\pi i}\bm b_{n,n+i-1}^T \\.
            2\pi i\bm a_{n-1,n-i}^T \bm \gamma_{n-1} & \bm a_{n,n-i}^T \\
         \end{MAT} \right),\qquad i\geq 0.
\end{gather}
\end{Corollary}

\begin{Remark}
One of the consequences of \eqref{symmetryForYnew} is that $\bm a_{n,j}$ and $\bm b_{n,j}$
appearing in~\eqref{Yni},~\eqref{Ynig} belong to $\RR^{N\times N}$.
\end{Remark}

The explicit expressions for coef\/f\/icients of the asymptotic expansion above, combined with the obvious fact that $\bm Y^n  (\bm Y^n)^{-1}=\bm I_{2N}$,  yield in a straightforward way the following identities:

\begin{Proposition}\label{lemmaab}
The coefficients of the monic MOPRL $\widehat{\bm P}_n$ and the coefficients of the Cauchy transform $\mathcal{C}(\widehat{\bm P}_n\bm W)$ $($i.e.\ the sequences $\bm b_{n,k})$ satisfy the following relations
\begin{gather}\label{anbn1}
    \sum_{j=0}^m \bm a_{n,n-m+j}\bm b_{n-1,n+j-1}^*=\sum_{j=0}^m\bm b_{n,n+m-j-1}\bm a_{n-1,n-j}^*,\qquad m\geq1,
\end{gather}
and
\begin{gather}\label{anbn2}
\sum_{j=0}^m \bm a_{n,n-m+j}\bm b_{n,n+j-1}^*=\sum_{j=0}^m\bm b_{n,n+m-j-1}\bm a_{n,n-j}^*,\qquad m\geq1.
\end{gather}
\end{Proposition}

In particular, for $m=1$ we recover the second identity of Corollary~\ref{corobn}.
\begin{proof}
From the asymptotic behavior of $\bm Y^n$ and $(\bm Y^n)^{-1}$ and using $ \bm Y^{n} (\bm Y^{n})^{-1}=\bm I$ we get
\begin{gather*}
\sum_{j=0}^m \bm Y^n_{m-j}\widetilde{\bm Y}^n_j=\bm 0,\qquad m\geq1,
\end{gather*}
from which block entries $(1,1)$ and $(1,2)$ give \eqref{anbn1} and \eqref{anbn2}, respectively.
\end{proof}

We f\/inish this section with the matrix analogue of the Christof\/fel--Darboux (CD) formula for MOPRL, written in terms of the solution of the characterizing RHP. For an orthonormal family  $(\bm P_n)_n$ we consider the kernel
\begin{gather}\label{Kg}
   \bm  K_n(x,y)=\sum_{j=0}^{n-1}\bm P_j^*(y)\bm P_j(x),\qquad x,y\in \mathbb R.
\end{gather}
Then the CD formula (see \cite{MR1426899}) reads as
\begin{gather}\label{CD}
   \bm  K_n(x,y)=\frac{\bm P_{n-1}^*(y)\bm A_n\bm P_n(x)-\bm P_n^*(y)\bm A_n^*\bm P_{n-1}(x)}{x-y}.
\end{gather}
\begin{Proposition}
The matrix CD kernel \eqref{Kg}
satisfies, for $x, y \in \mathbb R$,
\begin{gather*}
    \bm K_n(x,y) =-\frac{1}{2\pi
   i(x-y)}(\bm Y_{21}^*(y)\bm Y_{11}(x)-\bm Y_{11}^*(y)\bm Y_{21}(x))\\
\phantom{\bm K_n(x,y)}{} =\frac{1}{2\pi
   i(x-y)}\begin{pmatrix}
    \bm  0 &\bm  I
   \end{pmatrix}\big(\bm Y^{n}\big)_{+}^{-1}(y)\big(\bm Y^{n}\big)_{+}(x)\begin{pmatrix}
    \bm  I \\
     \bm 0 \\
   \end{pmatrix}.%\label{KRH}
\end{gather*}
\end{Proposition}

The proof follows taking into account the Christof\/fel--Darboux formula
(\ref{CD}), the expression for $\bm A_n$ in (\ref{An}) and the
expression of the inverse $(\bm Y^n)^{-1}$ in (\ref{RHSi}). The last proposition for general weight matrices of rectangular size can be found in \cite{MR2610344}.

\subsection{Dif\/ference and dif\/ferential relations}\label{ttrrsec}

Both the MOPRL and the solution of the RHP given in Theorem \ref{thm} depend on two variables, the discrete $n$ and the continuous $z$. We can derive further properties by analyzing the variation of this solution with respect to either variable, which yields linear relations of the form
\begin{gather}\label{LaxP}
        \bm Y^{n+1}(z)  =   \bm  E_n(z) \bm Y^{n}(z),\qquad
        \D\frac{d}{dz} \bm Y^{n}(z)  =   \bm F_n(z) \bm Y^{n}(z).
      \end{gather}
The general methodology to get these equations may be traced back to the original work of Gelfand, Levitan and other authors. The idea is that if the jump matrix for a RHP is independent of a variable, then a variation with respect to that variable leads to an identity. The fact that the jump matrix in~\eqref{jump_cond} is independent of $n$ allows one to obtain immediately the f\/irst identity in~\eqref{LaxP}, which in turn is connected to the well-known  three-term recurrence relation for MOPRL.

\begin{Theorem}\label{LaxP1}
There is a unique sequence of matrix coefficients $\bm \alpha_n$ such that the solution of the RHP $\bm Y^n(z)$ satisfies the following first-order difference equation
\begin{gather}\label{TTRRs}
\bm Y^{n+1}(z)=\left(\begin{MAT}{c.c}
  z\bm I-\bm \alpha_n & \frac{1}{2\pi i}\bm \gamma_n^{-1} \\.
  -2\pi i\bm \gamma_n & \bm 0 \\
\end{MAT} \right) \bm Y^{n}(z).
\end{gather}
Consequently, the monic MOPRL $\widehat{\bm P}_{n}$ satisfy the following three-term recurrence relation:
\begin{gather}\label{ttrrm}
x\widehat{\bm P}_n(x)=\widehat{\bm P}_{n+1}(x)+\bm \alpha_n\widehat{\bm P}_n(x)+\bm \beta_n\widehat{\bm P}_{n-1}(x),\qquad
n\geq0.
\end{gather}
The coefficients $\bm \alpha_n$ and $\bm \beta_n$, as well as $\bm \gamma_n$ $($defined in~\eqref{def:gamma}$)$ can be expressed either in terms of the elements of the solution of the RHP or in terms of the coefficients $\bm a_{n,j}$  $($see \eqref{orthogonality}$)$ as follows:
\begin{gather}
\bm \alpha_n  = \big(\bm Y_1^{n}\big)_{11}-\big(\bm Y_1^{n+1}\big)_{11}= \bm a_{n,n-1} - \bm a_{n+1,n}, \qquad
  \bm \beta_n  = \big(\bm Y_1^{n}\big)_{12}\big(\bm Y_1^{n}\big)_{21}=\bm \gamma_n^{-1}\bm \gamma_{n-1},\nonumber \\
  \bm \gamma_{n}  = -\frac{1}{2\pi i}\big(\bm Y_1^{n+1}\big)_{21}=-\frac{1}{2\pi i}\big(\bm Y_1^{n}\big)_{12}^{-1} .\label{betanttrr}
\end{gather}
\end{Theorem}

\begin{proof}
The matrix-valued function $\bm R=\bm Y^{n+1}(\bm Y^{n})^{-1}$, analytic in
$\mathbb{C}\setminus[a, b]$, satisf\/ies $\bm R_+(x)=\bm R_-(x)$ for all $x\in(a\,, b\, )$. By the Morera's theorem, $\bm R$ is analytic in $\mathbb{C}\setminus\{a,b\}$. The behavior at the exceptional points gives that the singularities at $a$ and $b$ are removable, and $\bm R$ is an entire function. From \eqref{Yninf} and Liouville's theorem we conclude that
\begin{gather}\label{TTRR}
\bm Y^{n+1}(z)= \left(\begin{MAT}{c.c}
  z\bm I+\big(\bm Y_1^{n+1}\big)_{11}-\big(\bm Y_1^{n}\big)_{11} & -\big(\bm Y_1^{n}\big)_{12} \\.%[2mm]
  \big(\bm Y_1^{n+1}\big)_{21} & \bm 0 \\
 \end{MAT} \right) \bm Y^{n}(z).
\end{gather}
This proves \eqref{TTRRs} along with the expression for $\bm \alpha_n$ in terms of the solution of the RHP. The three term recurrence relation~\eqref{ttrrm} for the MOPRL is obtained by considering the $(1,1)$ block entry of \eqref{TTRR}. Finally, the rest of the expressions for the coef\/f\/icients is found using  Corollary~\ref{coro25} (formula~\eqref{Yni}).
\end{proof}

We turn now to the dependence of $\bm Y^{n}$ on the continuous variable $z$; notice however that the jump \eqref{jump_cond} along the real line does depend on this variable, so the application of the paradigm explained at the beginning of this subsection is not straightforward. At this stage we can only aspire to f\/ind a dif\/ferential relation like in~\eqref{LaxP}, but with a matrix $\bm F_n$ depending upon the entries of $\bm Y^{n}$ (or equivalently, upon the MOPRL themselves).
\begin{Theorem}\label{TeordiffY}
The solution of the RHP for $\bm Y^n(z)$ satisfies the following first-order matrix differential equation
\begin{gather}\label{Deq}
    \frac{d}{dz}\bm Y^n(z)= \left(\begin{MAT}{c.c}
  -\mathfrak{B}_n(z) & -\frac{1}{2\pi i}\bm \gamma_{n}^{-1} \, \mathfrak{A}_n(z) \\.%[1mm]
 2\pi i\, \mathfrak{A}_{n-1}(z)\bm \gamma_{n-1} & \mathfrak{B}^*_n(z) \\
 \end{MAT} \right) \bm Y^{n}(z),
\end{gather}
where
\begin{gather}\label{Anz}
\mathfrak{A}_n(z)=-\bm \gamma_n\left(\D\int_a^b\D\frac{\widehat{\bm P}_n(t)\bm W'(t)\widehat{\bm P}^*_{n}(t)}{t-z}dt-\D\frac{\widehat{\bm P}_n(t)\bm W(t)\widehat{\bm P}^*_{n}(t)}{t-z}\bigg|_{t=a}^{t=b}\right)
\end{gather}
and
\begin{gather}\label{Bnz}
\mathfrak{B}_n(z)=-\left(\D\int_a^b\D\frac{\widehat{\bm P}_n(t)\bm W'(t)\widehat{\bm P}^*_{n-1}(t)}{t-z}dt-\D\frac{\widehat{\bm P}_n(t)\bm W(t)\widehat{\bm P}^*_{n-1}(t)}{t-z}\bigg|_{t=a}^{t=b}\right)\bm \gamma_{n-1},
\end{gather}
provided that all the functions in the right-hand sides of \eqref{Anz} and \eqref{Bnz} exist for $z\in\mathbb{C}\setminus[a, b]$. In particular,
\begin{gather}\label{RelA}
    \bm \gamma_{n}  \mathfrak{A}_n^*(z)=\mathfrak{A}_n(z)  \bm \gamma_{n}.
\end{gather}
\end{Theorem}

\begin{Remark}
Formulas (\ref{Anz}) and (\ref{Bnz}) were introduced in the scalar case in \cite{MR1616931} and in the matrix case (for a normalized family) in~\cite{MR2209517}. In general, they bear a formal character. In particular, the boundary terms in~(\ref{Anz}) and~(\ref{Bnz}) may be undef\/ined, depending on the behavior of~$\bm W$ at the endpoints (see (\ref{BouW})). As it will be shown in the next section, these formulas simplify considerably after considering a special transformation of the RHP for~$\bm Y^n$. Finally observe that the coef\/f\/icients~(\ref{Anz}) and~(\ref{Bnz}) are unique since $\bm Y^n$ is unique.
\end{Remark}

\begin{proof}
Consider the matrix-valued function $\bm R^n(z)=\left[\frac{d}{dz}\bm Y^n(z)\right](\bm Y^n)^{-1}(z)$. From the explicit expressions (\ref{RHS1}) and (\ref{RHSi}) we have that the block entries of $\bm R^n$ are
\begin{gather*}
(\bm R^n(z))_{11}  = -2\pi i\big(\widehat{\bm P}_n'(z)\mathcal{C}(\bm W\widehat{\bm P}_{n-1}^*)(z)-\mathcal{C}(\widehat{\bm P}_n\bm W)'(z)\widehat{\bm P}_{n-1}^*(z)\big)\bm \gamma_{n-1}, \\
(\bm R^n(z))_{12}  = -\widehat{\bm P}_n'(z)\mathcal{C}(\bm W\widehat{\bm P}_{n}^*)(z)+\mathcal{C}(\widehat{\bm P}_n\bm W)'(z)\widehat{\bm P}_{n}^*(z), \\
(\bm R^n(z))_{21}  = (2\pi i)^2\bm \gamma_{n-1} \big(\widehat{\bm P}_{n-1}'(z)\mathcal{C}(\bm W\widehat{\bm P}_{n-1}^*)(z)-\mathcal{C}(\widehat{\bm P}_{n-1}\bm W)'(z)\widehat{\bm P}_{n-1}^*(z)\big)\bm \gamma_{n-1},\\
(\bm R^n(z))_{22}  = -2\pi i\bm \gamma_{n-1}\big(-\widehat{\bm P}_{n-1}'(z)\mathcal{C}(\bm W\widehat{\bm P}_{n}^*)(z)+\mathcal{C}(\widehat{\bm P}_{n-1}\bm W)'(z)\widehat{\bm P}_{n}^*(z)\big).
\end{gather*}
In particular, $-(\bm R^n(z))_{11}=(\bm R^n(z))_{22}^*$ and $-2\pi i\bm\gamma_n(\bm R^n(z))_{12}=\frac{1}{2\pi i}(\bm R^{n+1}(z))_{21}\bm\gamma_n^{-1}$.

We will now derive formulas (\ref{Anz}) and (\ref{Bnz}). For that purpose we need the following technical observations:
\begin{Lemma} \label{lemma:2.3}\qquad
\begin{enumerate}\itemsep=0pt
\item[$(i)$] For every $\bm P\in\mathbb{P}_n$,
    \begin{gather*}\label{1p}
       \bm P(z)\mathcal{C}(\bm W\widehat{\bm P}_n^*)(z)=\mathcal{C}(\bm P\bm W\widehat{\bm P}_n^*)(z), \qquad
       \mathcal{C}(\widehat{\bm P}_n\bm W)(z)\bm P(z)=\mathcal{C}(\widehat{\bm P}_n\bm W\bm P)(z).
\end{gather*}
\item[$(ii)$]  For any differentiable and integrable matrix function $\bm F$ bounded at $z=a,b$,
      \begin{gather*}\label{2p}
\frac{d}{dz}  \mathcal{C}(\bm F)(z)=\mathcal{C}(\bm F')-\frac{1}{2\pi i}\frac{\bm F(t)}{t-z}\bigg|_{t=a}^{t=b}, \qquad z\in \mathbb{C}\setminus [a,b].
\end{gather*}
\end{enumerate}
\end{Lemma}

Identities in $(i)$ follow by adding and subtracting $\mathcal{C}(\bm P\bm W\widehat{\bm P}_n^*)(z)$ and using the orthogonality of the MOPRL, and $(ii)$ is obtained by integration by parts and using that $\frac{d}{dz}(1/(t-z))=-\frac{d}{dt}(1/(t-z))$.

We return to the proof of Theorem~\ref{TeordiffY}, showing in detail how to obtain $\mathfrak{B}_n(z)$. For $\mathfrak{A}_n(z)$ the computations are similar and will be omitted for the sake of brevity.

Using $(i)$ from the previous lemma and the product rule of dif\/ferentiation we get
\begin{gather*}
(\bm R^n(z))_{11}=-2\pi i\big(\mathcal{C}(\widehat{\bm P}_n'\bm W\widehat{\bm P}_{n-1}^*)(z)-\mathcal{C}(\widehat{\bm P}_n\bm W\widehat{\bm P}_{n-1}^*)'(z)+\mathcal{C}(\widehat{\bm P}_n\bm W(\widehat{\bm P}_{n-1}^*)')(z)\big)\bm \gamma_{n-1}.
\end{gather*}
Applying now $(ii)$ to the second term and canceling we get
\begin{gather*}
(\bm R^n(z))_{11}=-2\pi i\left(-\mathcal{C}(\widehat{\bm P}_n\bm W'\widehat{\bm P}_{n-1}^*)(z)+\frac{1}{2\pi i}\frac{\widehat{\bm P}_n(t)\bm W(t)\widehat{\bm P}_{n-1}^*(t)}{t-z}\bigg|_{t=a}^{t=b}\right)\bm \gamma_{n-1},
\end{gather*}
which yields~\eqref{Bnz}.
\end{proof}

Let us discuss now two main consequences of the existence of the f\/irst order matrix-valued dif\/ferential equation \eqref{Deq}.

One one hand, equations \eqref{LaxP}, known as the \emph{Lax pair} for the RHP, are clearly overdetermined, so compatibility conditions (via  cross-dif\/ferentiation of both equations) yield
\begin{gather}\label{Comp}
 \bm E_n'(z)+ \bm E_n(z) \bm F_n(z)= \bm F_{n+1}(z) \bm E_n(z),
\end{gather}
also known as \emph{string equations}. They allow us, for instance, to recover the sequence $\{\bm F_n\}$ from the sequence $\{\bm E_n\}$ and the initial value $\bm F_1$. These ideas relate also the RHP to nonlinear problems such as the nonlinear Schr\"{o}dinger equation (NLS), see~\cite{MR1989226}, and the Toda f\/low, see~\cite{MR0408647}. In our situation, the compatibility conditions~\eqref{Comp} considered entry-wise imply the following:

\begin{Proposition}\label{Stringcoro}
The recurrence relations
\begin{gather}\label{String11}
   \bm I+\mathfrak{B}_{n+1}(z)(z\bm I-\bm \alpha_n)-(z\bm I-\bm \alpha_n)\mathfrak{B}_n(z)=\mathfrak{A}_{n+1}^*(z)\bm \beta_{n+1}-\bm \beta_n\mathfrak{A}_{n-1}^*(z)
\end{gather}
and
\begin{gather}\label{String21}
    \mathfrak{B}_{n+1}(z)+\bm \gamma_{n}^{-1}\mathfrak{B}_{n}^*(z)\bm \gamma_{n}=(z\bm I-\bm \alpha_n)\mathfrak{A}_n^*(z)
\end{gather}
hold for every $n\geq0$, where $\bm\alpha_n$ and $\bm\beta_n$ are the coefficients of the three-term recurrence rela\-tion~\eqref{ttrrm}.
\end{Proposition}

\begin{proof}
Block entries $(1,1)$ and $(1,2)$ of (\ref{Comp}) give (\ref{String11}) and (\ref{String21}), respectively. Block entry $(2,1)$ of (\ref{Comp}) is equivalent to (\ref{String21}) using relation (\ref{RelA}) and $\bm \alpha_n=\bm\gamma_{n}^{-1}\bm \alpha_n^*\bm \gamma_{n}$ (which is a~consequence of substituting \eqref{Yni} for $i=2$ in \eqref{TTRRs} as $z\rightarrow\infty$, along with \eqref{betanttrr}), while block entry~$(2,2)$ is equivalent to~(\ref{String11}) as well.
\end{proof}

Another consequence of Theorem~\ref{TeordiffY} is the existence of the \emph{ladder operators} for MOPRL:

\begin{Corollary}\label{corolo}
The monic MOPRL $(\widehat{\bm P}_n)_n$ satisfy the following difference-differential  relations $($lowering and raising operators, respectively$)$:
\begin{gather}\label{lowop1}
    \widehat{\bm P}_n'(z)=-\mathfrak{B}_n(z)\widehat{\bm P}_n(z)+\mathfrak{A}_n^*(z)\bm \beta_n\widehat{\bm P}_{n-1}(z)
\end{gather}
and
\begin{gather}\label{raiop1}
    \widehat{\bm P}_n'(z)=[\mathfrak{A}_{n}^*(z)(z\bm I-\bm\alpha_n)-\mathfrak{B}_{n}(z)]\widehat{\bm P}_{n}(z)-\mathfrak{A}_n^*(z)\widehat{\bm P}_{n+1}(z).
\end{gather}
\end{Corollary}

\begin{proof}
Block entry $(1,1)$ of (\ref{Deq}) gives the lowering operator (\ref{lowop1}), using the expression for~$\bm \beta_n$ in Theorem~\ref{LaxP1}, while block entry $(2,1)$ of~(\ref{Deq}) gives the raising operator~(\ref{raiop1}) using~\eqref{RelA} and~\eqref{String21}.
\end{proof}

Proposition \ref{Stringcoro} and Corollary \ref{corolo} were already known for the normalized MOPRL in \cite{MR2209517}; the RH approach gives an alternative proof of these results. The ladder operators are the basic dif\/ferential relations for MOPRL. It is well-known that they can be combined to build a second-order dif\/ferential equation satisf\/ied by the polynomials. This fact was mentioned in~\cite{MR2209517}, but no explicit expression of the dif\/ferential equation was given; here we present it for completeness:
\begin{Corollary}\label{sodeq1}
The polynomials $(\widehat{\bm P}_n)_n$ satisfy the following second-order differential equation
\begin{gather}\label{seqex}
    \widehat{\bm P}_n''(z)+\mathfrak{M}_n(z)\widehat{\bm P}_n'(z)+\mathfrak{N}_n(z)\widehat{\bm P}_n(z)=\bm 0,
\end{gather}
where
\begin{gather*}\label{Mnz}
   \mathfrak{M}_n(z)=-(\mathfrak{A}_n^{*}(z))'\mathfrak{A}_n^{-*}(z)+\mathfrak{B}_n(z)-\mathfrak{A}_n^{*}(z)(z\bm I-\bm\alpha_n)+\mathfrak{A}_n^{*}(z)\mathfrak{B}_{n+1}(z)\mathfrak{A}_n^{-*}(z)
\end{gather*}
and
\begin{gather*}\label{Nnz}
   \mathfrak{N}_n(z)=\mathfrak{M}_n(z)\mathfrak{B}_n(z)-\mathfrak{B}_n^2(z)+\mathfrak{B}_n'(z)+\mathfrak{A}_n^{*}(z)\bm\beta_n\mathfrak{A}_{n-1}^{*}(z),
\end{gather*}
provided that the inverse of $\mathfrak{A}_n^{*}(z)$ exists for $z\in\mathbb{C}\setminus[a, b]$. Here and below, the superscript $-*$ denotes the conjugate transpose of the inverse.
\end{Corollary}

\begin{proof}
Dif\/ferentiate the lowering operator (\ref{lowop1}) and substitute the raising operator (\ref{raiop1}) evaluated at $n-1$.
\end{proof}

\begin{Remark}
Note that we can easily obtain another second-order dif\/ferential operator satisf\/ied by the orthogonal polynomials reversing the ladder operators. This new dif\/ferential equation needs not in principle be the same as (\ref{seqex}). Nevertheless, it is straightforward to see that both equations are equivalent using the compatibility conditions (\ref{String11}) and (\ref{String21}).
\end{Remark}

Although Theorem \ref{TeordiffY} gives explicit expressions of $\mathfrak{A}_n(z)$ and $\mathfrak{B}_n(z)$, these coef\/f\/icients are usually dif\/f\/icult to calculate. They are not even def\/ined on the interval $[a, b]$. In the next section we will introduce some additional assumptions that simplify the dif\/ferential equation considerably.

\section[Transformation of the RHP when ${\rm supp}(W)=\mathbb{R}$]{Transformation of the RHP when $\boldsymbol{{\rm supp}(W)=\mathbb{R}}$}\label{secRH}

As it was pointed out in the previous section, the independence of the jump in \eqref{jump_cond} with respect to the discrete variable $n$ allows one to obtain the three-term recurrence relation in a straightforward way. Under additional assumptions on the jumps we can perform a transformation of the RHP in such a way that a new jump is independent of the continuous variable $z$. This will have consequences on the resulting dif\/ferential relations.

In the matrix case there is some extra freedom absent in the scalar situation, that gives us a~whole new family of dif\/ferential relations, and consequently, a whole \emph{class} of ladder operators, some of them of the $0$-th order, something that is not possible in the scalar case. The combination of the ladder operators will give rise to a \emph{class} of second-order dif\/ferential equations, some of them of order less than $2$.
%of first and 0-th order.

For simplicity, we will focus here on the case when $[a,b]=\mathbb{R}$, assuming along this section that~$\bm W$ is smooth and positive def\/inite on the whole real axis $\RR$. Some ideas about how to handle the case of f\/inite endpoints of the support of $\bm W$ are brief\/ly explained in the f\/inal Section~\ref{sec:finalremarks}.

\subsection{The transformation}\label{sec:transf}

As we have seen in Section~\ref{ttrrsec},  the recurrence relation (f\/irst identity in~\eqref{LaxP}) is a general fact intrinsic to the orthogonality of the polynomials. Simple dif\/ferential relations are much more demanding, and we must impose at this stage further conditions on the orthogonality weight~$\bm W$.

Our immediate goal is to obtain an invertible transformation $\bm Y^{n} \rightarrow \bm X^{n}$ such that $\bm X^n$ has a~constant jump across $\RR$. We will consider
\begin{gather*}
\bm X^{n}(z) :=\bm Y^{n}(z) \bm V(z),
\end{gather*}
where $\bm V$ is a matrix-valued function, analytic in $\CC\setminus \RR$ and continuous up to $\mathbb{R}$, and invertible for all $z\in \CC$. Then the jump matrix for $\bm X^n$  on $\RR$ is
\begin{gather} \label{motivation1}
\bm V^{-1} \left(\begin{MAT}{c.c}
                             \bm   I & \bm W \\.
                            \bm    0 &\bm  I \\
                             \end{MAT} \right) \bm V.
\end{gather}
Observe that this kind of transformations does not af\/fect the f\/irst dif\/ference equation in the Lax pair~(\ref{LaxP}), but constant jumps will allow us to use the strategy of variation of the problem along $z$. Going down the path of simplif\/ication, we can try a block-diagonal matrix
\begin{gather*}
\bm V(z)= \left(\begin{MAT}{c.c}
     \bm T(z)& \bm 0 \\.
      \bm 0 & \bm T^{-*}(z)\\
     \end{MAT} \right) ,
\end{gather*}
where $\bm T$ is an invertible $N\times N$ matrix function. Then \eqref{motivation1} boils down to
\begin{gather*}
 \left(\begin{MAT}{c.c}
                             \bm   I & \bm T^{-1} \bm W \bm T^{-*} \\.
                            \bm    0 &\bm  I \\
                             \end{MAT} \right) .
\end{gather*}

These preliminary considerations motivate to consider a factorization of the weight in the form
\begin{gather}\label{Ww}
\bm W(x)= {\bm T}(x) {\bm T}^*(x),\qquad x\in\mathbb{R},
\end{gather}
where $\bm T$ is a smooth matrix-valued function on $\RR$. Under our assumptions the existence of such a $\bm T$ is guaranteed, but not its uniqueness: performing its $RQ$ decomposition it can be written~as
\begin{gather}\label{WwT}
 \bm T(x)=\widehat{\bm T}(x)\bm S(x) ,\qquad   x\in\mathbb{R} ,
\end{gather}
where $\widehat{\bm T}(x)$ is an upper triangular matrix and $\bm S(x)$ is an arbitrary smooth and unitary matrix (for each $x\in \mathbb R$). This representation can go further taking into account that any unitary matrix can be written as $\bm S(x)=e^{\bm Q (x)}$, where $\bm Q(x)$ is an skew-Hermitian matrix function. Additionally, since any skew-Hermitian matrix is normal, it has a factorization $\bm Q(x)=i\bm U(x)\bm D(x)\bm U^*(x)$, where $\bm U$ is again unitary and $\bm D$ a diagonal matrix with real entries.
Therefore $\bm T$ from \eqref{WwT} can be written as
\begin{gather*}
 \bm T(x)=\widehat{\bm T}(x)\bm U(x)e^{i\bm D(x)}\bm U^*(x).
\end{gather*}

We narrow the choice of $\bm T$ (and $\bm W$) by imposing the additional constraint that there exists a matrix polynomial $\bm G$ such that in a neighborhood of the origin,
\begin{gather}\label{ode:G}
\bm T'(z)= {\bm G}(z) {\bm T} (z).
\end{gather}
As a consequence, (see the discussion in \cite[Chapter~9]{MR0499382} or \cite[Chapter~1]{Sibuya:1990}), $\bm T$ has an analytic continuation to the whole plane as an entire and invertible matrix-valued function, and{\samepage
\begin{gather*}
\label{analyticExtW}
\bm W (z)=\bm T(z) \bm T^{*}(z)
\end{gather*}
provides an analytic extension of $\bm W$ to the whole plane.}

Following our previous discussion, we def\/ine the entire matrix-valued function
\begin{gather*}
\bm V(z)= \left(\begin{MAT}{c.c}
     \bm T(z)& \bm 0 \\.
      \bm 0 & \bm T^{-*}(z)\\
    \end{MAT} \right)  .
\end{gather*}
Observe that $\bm V$ is invertible for all $z\in \CC$. By dif\/ferentiating the identity $\bm T^{-1} \bm T =\bm I$ and using~\eqref{ode:G} we conclude that $\left( \bm T^{-1}\right)'=-  \bm T^{-1} \bm G$, so that
\begin{gather}
\label{diffEqforV}
\bm V'(z)= \left(\begin{MAT}{c.c}
        \bm G(z) & \bm 0\\.
        \bm 0 & -\bm G^*(z)\\
       \end{MAT} \right) \bm V(z).
\end{gather}

If $\bm Y^n(z)$ is the solution of the original RHP, let us def\/ine
\begin{gather}\label{X}
    \bm X^n(z)=\bm Y^n(z)  \bm V(z).
\end{gather}
Then straightforward computations show that $\bm X^n(z)$ is analytic in
$\mathbb{C}\setminus\mathbb{R}$, satisf\/ies the jump condition
\begin{gather*} \label{constantJumps}
\bm X^n_+(x)=\bm X^n_-(x) \left(\begin{MAT}{c.c}
  \bm I_N & \bm I_N \\.
  \bm 0 & \bm I_N \\
\end{MAT} \right) , \qquad x\in \RR,
\end{gather*}
and for $m\in \NN$ (see \eqref{Yninf}),
\begin{gather}\label{Xinff}
\bm X^n(z)=\left(\bm I_{2N}+\sum_{i=1}^m \dfrac{\bm Y_i^n}{z^i}+\mathcal{O} (1/z^{m+1})\right) \left(\begin{MAT}{c.c}
  z^n\, \bm I_N & \bm 0 \\.
  \bm 0 & z^{-n}\, \bm I_N \\
\end{MAT} \right)  \bm V(z), \qquad z\to \infty.
    \end{gather}
Since $\bm X^n(z)$ is invertible in $\CC \setminus \RR$,  we can consider the matrix function $\bm F_n(z){=}\big[\frac{d}{dz}\bm X^n(z)\big]\!\bm X^n(z)^{-1}.\!$ Again, $\bm F_n(z)$ is analytic in $\mathbb{C}\setminus\mathbb{R}$ and on the real line,
$(\bm F_n)_+(x)=(\bm F_n)_-(x)$, which implies that it is an entire matrix function.
From \eqref{Yninf} and (\ref{Xinff}), for $z\to \infty$,
\begin{gather*}\label{dXinf}
 \left[ \frac{d}{dz}\bm Y^n(z) \right] \left[  \bm Y^n(z) \right]^{-1}=\mathcal O\left( \frac{1}{z}\right),
\end{gather*}
and combining it with \eqref{Yninf}, \eqref{Yninfi}, \eqref{diffEqforV} and \eqref{Xinff}, we get for $z\rightarrow\infty$,
\begin{gather}
\bm F_n(z)  =\left[\frac{d}{dz}\bm X^n(z)\right]\bm X^n(z)^{-1}\nonumber \\
\phantom{\bm F_n(z)}{}  =\left(\bm I_{2N}+\sum_{i=1}^m \dfrac{\bm Y_i^n}{z^i}+\mathcal{O} \big(1/z^{m+1}\big)\right) \left(\begin{MAT}{c.c}
        \bm G(z) & \bm 0\\.
        \bm 0 & -\bm G^*(z)\\
      \end{MAT} \right)  \nonumber\\
\phantom{\bm F_n(z)=}{} \times \left(\bm I_{2N}+ \sum_{i=1}^m \dfrac{\widetilde{\bm Y}_i^n}{z^i}+\mathcal{O} \big(1/z^{m+1}\big)\right)+\mathcal{O}(1/z).\label{Xinf}
\end{gather}
By Liouville's theorem, the right hand side in \eqref{Xinf} will coincide with its polynomial part in the expansion at inf\/inity, and its degree is no greater than the degree of $\bm G$. To be more precise, if we assume that the degree of $\bm G$ is $m\in \NN$ and denote
\begin{gather*}
\bm G(z)=\sum_{j=0}^m\bm M_jz^j, \qquad \widetilde{\bm M}_j= \left(\begin{MAT}{c.c}
       \bm M_j & \bm 0\\.
       \bm 0 & -\bm M_j^*\\
     \end{MAT} \right)  ,
\end{gather*}
then after dropping the negative powers of $z$ in
\begin{gather*}
\bm F_n(z)   =\left(\bm I_{2N}+\sum_{i=1}^m \dfrac{\bm Y_i^n}{z^i}+\mathcal{O} \big(1/z^{m+1}\big)\right)\\
 \hphantom{\bm F_n(z)   =}{}\times \left( \sum_{j=0}^m\widetilde{\bm M}_jz^j \right)
 \left(\bm I_{2N}+ \sum_{k=1}^m \dfrac{\widetilde{\bm Y}_k^n}{z^k}+\mathcal{O} \big(1/z^{m+1}\big)\right)+\mathcal{O}(1/z)
\end{gather*}
we obtain that
\begin{gather}\label{Rpolyn}
\bm F_n(z)= \sum_{k=0}^m \left( \sum_{j=k}^m \sum_{i=0}^{j-k} \bm Y_i^n \widetilde{\bm M}_j \widetilde{\bm Y}_{j-i-k}^n \right) z^k, \qquad \text{with} \quad \bm Y_0^n=\widetilde{\bm Y}_{0}^n=\bm I_{2N}.
   \end{gather}

For instance, if $m=0$,
\begin{gather}\label{Rm0}
    \bm F_n(z)=\widetilde{\bm M}_0,
\end{gather}
for $m=1$,
\begin{gather}\label{Rm1}
    \bm F_n(z)=\widetilde{\bm M}_1z+\widetilde{\bm M}_0+\bm Y^n_1\widetilde{\bm M}_1+\widetilde{\bm M}_1\widetilde{\bm Y}^n_1,
\end{gather}
and for $m=2$,
\begin{gather}
    \bm F_n(z)=   \widetilde{\bm M}_2z^2+(\widetilde{\bm M}_1+\bm Y^n_1\widetilde{\bm M}_2+\widetilde{\bm M}_2\widetilde{\bm Y}^n_1)z+\widetilde{\bm M}_0+\bm Y^n_1\widetilde{\bm M}_1+\widetilde{\bm M}_1\widetilde{\bm Y}^n_1\nonumber\\
\phantom{\bm F_n(z)=}{}  +\bm Y^n_2\widetilde{\bm M}_2+\bm Y^n_1\widetilde{\bm M}_2\widetilde{\bm Y}^n_1+\widetilde{\bm M}_2\widetilde{\bm Y}^n_2.\label{Rm2}
\end{gather}

Finally, using the explicit expressions (\ref{Yni}) and (\ref{Ynig}) in \eqref{Rpolyn}, we arrive at the following

\begin{Theorem}\label{thmXnp}
Under assumptions \eqref{Ww} and \eqref{ode:G}, with $\bm G(z)=\sum\limits_{j=0}^m\bm M_jz^j\in \mathbb P_m$,
the matrix function $\bm X^n$ defined in \eqref{X} satisfies the following first-order differential equation with polynomial coefficients:
\begin{gather}\label{DeqX}
    \frac{d}{dz}\bm X^n(z)= \bm F_n(z;\bm G) \bm X^{n}(z),
\end{gather}
where
\begin{gather}\label{DefX}
\bm F_n(z;\bm G) = \left(\begin{MAT}{c.c}
  -\mathcal{B}_n(z; \bm G) & -\frac{1}{2\pi i}\bm \gamma_{n}^{-1} \mathcal{A}_n(z; \bm G) \\.
 2\pi i\, \mathcal{A}_{n-1}(z; \bm G)\bm \gamma_{n-1} & \mathcal{B}^*_n(z; \bm G) \\
\end{MAT} \right) ,
\end{gather}
$\mathcal{A}_n$ and $\mathcal{B}_n$ are matrix polynomials,
\begin{gather}\label{AnzX}
\mathcal{A}_n(z; \bm G)  =-\bm\gamma_n\left(\sum_{j=0}^{m-1}\bm b_{n,n+m-j-1}\bm \Delta_{j,n}^*(z)+\bm \Delta_{j,n}(z)\bm b_{n,n+m-j-1}^*\right), \\
\label{BnzX}
\mathcal{B}_n(z; \bm G)  =-\left(\sum_{j=0}^{m}\bm\Delta_{j,n}(z)\bm b_{n-1,n+m-j-1}^*+\bm b_{n,n+m-j-2}\bm \Delta_{j,n-1}^*(z)\right)\bm\gamma_{n-1},
\end{gather}
and the coefficients $\bm \Delta_{j,n}(z)$ are given by
\begin{gather*}
\bm \Delta_{j,n}(z)=\sum_{k=0}^j\widehat{\bm P}_{n,k}(z)\bm M_{m-j+k},\qquad\widehat{\bm P}_{n,k}(z)=z^k\bm I+\bm a_{n,n-1}z^{k-1}+\cdots+\bm a_{n,n-k}.
\end{gather*}
Moreover, $\bm F_n(\cdot;\bm G) $ is linear in $\bm G$:
\begin{gather}\label{linearityF}
\bm F_n(z;\bm G_1+\bm G_2) =\bm F_n(z;\bm G_1)+ \bm F_n(z;\bm G_2).
\end{gather}
\end{Theorem}

\begin{Remark}
The coef\/f\/icients $\bm b_{n,k}$ were introduced in \eqref{bnk} and discussed in Lemma~\ref{lemmabn}. It is worth observing the similarity of this result with Theorem \ref{TeordiffY}; in the present situation we can compute the dif\/ferential equation for $\bm X^n$ directly in terms of the coef\/f\/icients of the MOPRL~$\widehat{\bm P}_n$ without considering integrals or studying the behavior at the endpoints. Finally, we have once again that
\begin{gather}\label{Angn}
\bm \gamma_{n}\mathcal{A}_n^*(z; \bm G)=\mathcal{A}_n(z; \bm G)\bm \gamma_{n}.
\end{gather}
\end{Remark}

Again, for lowest degrees $m$ in Theorem \ref{thmXnp} we can use (\ref{Rm0}), (\ref{Rm1}), (\ref{Rm2}), Lemma \ref{lemmabn} and Corollary \ref{corobn} in order to write the coef\/f\/icient matrix \eqref{DefX} of the dif\/ferential equation explicitly. For instance, for $m=0$,
\begin{gather*}\label{Rm0x}
\bm F_n(z;\bm G)= \left(\begin{MAT}{c.c}
                         \bm M_0 & \bm 0 \\.
                      \bm 0  &  -\bm M_0^*\\
                     \end{MAT} \right)  ;
\end{gather*}
for $m=1$,
\begin{gather}\label{Rm1x}
        \bm F_n(z;\bm G)= \!\left(\begin{MAT}{c.c}
                        \bm G(z) +\bm a_{n,n-1}\bm M_1-\bm M_1\bm a_{n,n-1}  & \frac{1}{2\pi i}(\bm \gamma_{n}^{-1}\bm M_1^*+\bm M_1\bm \gamma_{n}^{-1}) \\.
                       -2\pi i(\bm \gamma_{n-1}\bm M_1+\bm M_1^*\bm \gamma_{n-1})   &  -\bm G^*(z)+\bm a_{n,n-1}^*\bm M_1^*-\bm M_1^*\bm a_{n,n-1}^*\\
                      \end{MAT} \right)\!
 ,\!\!\!
\end{gather}
and for $m=2$, with the notation \eqref{DefX},
\begin{gather*}
\mathcal{B}_n(z; \bm G)= -\bm G(z)  -( \bm a_{n,n-1}\bm M_2-\bm M_2\bm a_{n,n-1})z +\bm M_1\bm a_{n,n-1}-\bm a_{n,n-1}\bm M_1-\bm a_{n,n-2}\bm M_2\nonumber\\
 \phantom{\mathcal{B}_n(z; \bm G)=}{}
 -\bm M_2(\bm a_{n+1,n}\bm a_{n,n-1}+\bm a_{n+1,n-1})+\bm a_{n,n-1}\bm M_2\bm a_{n,n-1}-\bm\gamma_{n}^{-1}\bm M_2^*\bm\gamma_{n-1},
\nonumber\\
\mathcal{A}_n(z; \bm G)= -\bm M_2^*z-\bm M_1^*+\bm a_{n+1,n}^*\bm M_2^*-\bm M_2^*\bm a_{n,n-1}^* \\
\phantom{\mathcal{A}_n(z; \bm G)=}{}
 -\bm\gamma_{n}(\bm M_2z+\bm M_1-\bm M_2\bm a_{n+1,n}+\bm a_{n,n-1}\bm M_2)\bm\gamma_{n}^{-1}.%\label{Rm2x}
\end{gather*}

Observe from Theorem \ref{thmXnp} that $\mathcal{B}_n(z; \bm G)$ is a matrix polynomial of degree at most~$m$ and~$\mathcal{A}_n(z; \bm G)$ is a matrix polynomial of degree at most $m-1$.

Finally, it is important to emphasize that in many situations we can exploit the implications of the freedom in the factorization~\eqref{WwT} on the dif\/ferential equation \eqref{DeqX}. This freedom, as we will see in Proposition~\ref{propsc}, appears only in the matrix setting.

\begin{Proposition}\label{thm:nonunique}
If under assumptions \eqref{Ww} and \eqref{ode:G}, with $\bm G$ a polynomial, there exists a~non-trivial matrix-valued function $\bm S$, non-singular on $\CC$, smooth and unitary on $\RR$, such that
\begin{gather}
\label{def:H}
\bm H(z)=\bm T(z)  \bm S'(z) \bm S^* (z) \bm T^{-1}(z)
\end{gather}
is also a polynomial, then $\widetilde{\bm T}=\bm T \bm S$ satisfies
\begin{gather*}
\bm W(x)= \widetilde{\bm T}(x) \widetilde{\bm T}^*(x),\quad x\in\mathbb{R}, \qquad \widetilde{\bm T}'(z)= \widetilde{\bm G}(z) \widetilde{\bm T} (z), \quad z\in \CC,
\end{gather*}
with $\widetilde{\bm G}(z)= \bm G(z)+ \bm H(z)$. Moreover, the matrix $\bm X^n$ defined in \eqref{X}, satisfies along with \eqref{DeqX}--\eqref{BnzX} the following relation:
\begin{gather}\label{DeqXbis}
    \frac{d}{dz}{\bm X}^n(z) =   \left( \bm F_n(z;\bm G)+\bm F_n(z;\bm H)\right) {\bm X}^n(z)  - {\bm X}^n(z)  \left(\begin{MAT}{c.c}
   \bm{\chi}(z) & \bm 0 \\.
      \bm 0 & -\bm{\chi}^*(z)\\
     \end{MAT} \right) ,
\end{gather}
with
$ \bm{\chi}(z)=      \bm S'(z)  \bm S^* (z)   $.
\end{Proposition}

\begin{Remark}
Equations  \eqref{DeqX} and  \eqref{DeqXbis} are not necessarily trivially related, and in principle we could combine them in order to establish new relations for $\bm X^n$, and thus, for $\bm Y^n$.
\end{Remark}

\begin{proof}
The key observation is the formula \eqref{linearityF}, so that the coef\/f\/icients $\mathcal A_n$ and $\mathcal B_n$ in \eqref{DeqX} depend linearly on $\bm G$:
\begin{gather}
\label{linearity}
\mathcal A_n(\cdot;  \bm G +\bm H)=\mathcal A_n(\cdot;{\bm G}) + \mathcal A_n(\cdot;{\bm H}), \qquad  \mathcal B_n(\cdot; \bm G +\bm H)=\mathcal B_n(\cdot;{\bm G}) + \mathcal B_n(\cdot;{\bm H}).
\end{gather}
\end{proof}

\subsection{Dif\/ferential properties}\label{sec:diffprop}

Taking into account the similarities between Theorems \ref{TeordiffY} and \ref{thmXnp}, we get immediately  the analogue of Proposition \ref{Stringcoro} (the compatibility conditions):
\begin{Proposition}\label{propcc}
Under assumptions \eqref{Ww} and \eqref{ode:G}, with $\bm G(z)$ a polynomial, the coefficients of the differential equation \eqref{DeqX} satisfy the following recurrence relations: for every $n\geq0$,
\begin{gather}\label{String11X}
   \bm I+\mathcal{B}_{n+1}(z;\bm G)(z\bm I-\bm\alpha_n)-(z\bm I-\bm\alpha_n)\mathcal{B}_n(z;\bm G)=\mathcal{A}_{n+1}^*(z;\bm G)\bm \beta_{n+1}-\bm\beta_n\mathcal{A}_{n-1}^*(z;\bm G)\!\!\!
\end{gather}
and
\begin{gather}\label{String21X}
    \mathcal{B}_{n+1}(z;\bm G)+\bm\gamma_{n}^{-1}\mathcal{B}_{n}^*(z;\bm G)\bm\gamma_{n}=(z\bm I-\bm\alpha_n)\mathcal{A}_n^*(z;\bm G),
\end{gather}
where $\bm\alpha_n$ and $\bm\beta_n$ are the coefficients of the three-term recurrence relation \eqref{ttrrm}.
\end{Proposition}

As before, the ladder operators (the most basic dif\/ferential properties for MOPRL) can be easily obtained by analyzing the f\/irst block column of $\bm X^n$  in the dif\/ferential equation \eqref{DeqX}:
\begin{Proposition}\label{propladdX}
Under assumptions \eqref{Ww} and \eqref{ode:G}, with $\bm G(z)$ a matrix polynomial, the monic MOPRL $(\widehat{\bm P}_n)_n$ satisfy the following difference-differential relations $($lowering and raising operators, respectively$)$:
\begin{gather}\label{lowop1X}
    \widehat{\bm P}_n'(z)+ \widehat{\bm P}_n(z) \bm G(z)= -\mathcal{B}_n(z;\bm G)\widehat{\bm P}_n(z)+\mathcal{A}_n^*(z;\bm G)\bm\beta_n\widehat{\bm P}_{n-1}(z)
\end{gather}
and
\begin{gather}\label{raiop1X}
    \widehat{\bm P}_n'(z)+\widehat{\bm P}_n(z) \bm G(z) =\big(\mathcal{A}_n^*(z;\bm G)(z\bm I-\bm\alpha_n)-\mathcal{B}_n(z;\bm G)\big)\widehat{\bm P}_{n}(z)-\mathcal{A}_n^*(z;\bm G)\widehat{\bm P}_{n+1}(z),
\end{gather}
where $\mathcal{A}_{n}$ and $\mathcal{B}_{n}$ are given by \eqref{AnzX} and \eqref{BnzX}, respectively.
\end{Proposition}

\begin{proof}
Block entry $(1,1)$ of (\ref{DeqX}) gives the lowering operator (\ref{lowop1X}), using \eqref{betanttrr} and \eqref{Angn}, while block entry $(2,1)$ of (\ref{DeqX}) gives the raising operator (\ref{raiop1X}) using \eqref{String21X} and \eqref{Angn}.
\end{proof}

\begin{Remark}
Comparing with the results of Corollary \ref{corolo}, notice that these ladder operators contain a term  with the MOPRL multiplied on the left.
\end{Remark}

In the situation described in Proposition \ref{thm:nonunique}, we can in principle exploit the non-uniqueness of the relations above in order to derive further relations for the MOPRL. The following result shows the possibility of the existence of 2-terms recurrence relations for the family $\widehat{\bm P}_n$, a phenomenon that has not been reported before in the theory of MOPRL.

\begin{Corollary} \label{coro:additional}
Under conditions of Proposition  {\rm \ref{thm:nonunique}},  the family of monic MOPRL $(\widehat{\bm P}_n)_n$ satisfies, along with \eqref{lowop1X}--\eqref{raiop1X}, the following  relations:
\begin{gather}\label{lowop2H}
\widehat{\bm P}_n(z)\bm H(z) =-\mathcal{B}_n(z; \bm H)\widehat{\bm P}_n(z)+\mathcal{A}_n^*(z;\bm H)\bm\beta_n\widehat{\bm P}_{n-1}(z)
\end{gather}
and
\begin{gather}\label{raiop2H}
\widehat{\bm P}_n(z)\bm H(z) =\big(\mathcal{A} _n^*(z; \bm H)(z\bm I-\bm\alpha_n)-\mathcal{B}_n(z; \bm H)\big) \widehat{\bm P}_{n}(z)-\mathcal{A} _n^*(z; \bm H)\widehat{\bm P}_{n+1}(z),
\end{gather}
where $\bm H$ is defined in \eqref{def:H}.

In particular, provided $\mathcal{A}_n(z;\bm G)$ is invertible for all $z\in \CC$,
\begin{gather}
\widehat{\bm P}_n(z)\bm H(z) +\mathcal{B}_n(z; \bm H)\widehat{\bm P}_n(z)-\mathcal{A}_n^*(z;\bm H) \mathcal{A}_n^{-*}(z;\bm G)\nonumber\\
\qquad{}\times
\big( \widehat{\bm P}_n'(z)+ \widehat{\bm P}_n(z) \bm G(z)+ \mathcal{B}_n(z;\bm G)\widehat{\bm P}_n(z) \big)=\bm 0.\label{ode1stOrder}
\end{gather}
\end{Corollary}

\begin{proof}
To get \eqref{lowop2H}, we subtract \eqref{DeqX} and \eqref{DeqXbis} and then evaluate the $(1,1)$ block entry. \eqref{raiop2H} is obtained from \eqref{lowop2H} and the three term recurrence relation. Finally, \eqref{ode1stOrder} is obtained by replacing \eqref{lowop1X} in \eqref{lowop2H}.
\end{proof}

Equations  \eqref{lowop2H} and \eqref{raiop2H} are known as the 0-th order ladder operators, while \eqref{ode1stOrder} is a~f\/irst order dif\/ferential relation for the MOPRL. In some situations they yield trivial identities; this is always true in the scalar case, as the following proposition shows:
\begin{Proposition}\label{propsc}
Assume that under the conditions of Proposition~{\rm \ref{thm:nonunique}},    $\bm\chi(z)= \bm S'(z) \bm S^* (z)=ip(z)\bm I$, where $p$ is a scalar polynomial of degree $m$. Then
\begin{gather*}
\mathcal{A}_n(z;\bm H)=\bm 0,\qquad\mbox{and}\qquad\mathcal{B} _n(z;\bm H)=-ip(z)\bm I.
\end{gather*}
\end{Proposition}

\begin{proof}
From \eqref{def:H}, $\bm H=\bm T \bm\chi \bm T^{-1}$, and by linearity of $\mathcal{A}_n$ and $\mathcal{B}_n$ in $\bm H$ (see~\eqref{linearity}), it is enough to prove the formulas above for monomials $z^k\bm I$. In this case, formulas~(\ref{AnzX}) and~(\ref{BnzX}) simplify considerably. Using Proposition~\ref{lemmaab} gives $\mathcal{A}_n(z;\bm H)=\bm 0$ and $\mathcal{B}_n(z;\bm H)=-iz^k\bm I$.
\end{proof}

In the scalar case ($N=1$) the only smooth and unitary function $s(x)$ on $\RR$ has the form $s(z)=e^{i p(x)}$, with $p$ real-valued, and
\begin{gather*}
s'(z)\overline{s(\bar{z})}=i p(z),
\end{gather*}
so that the assumption that $p$ is a polynomial brings us to the situation described in Proposition~\ref{propsc}. This explains why in the scalar setting we never get nontrivial 0-th order ladder operators. In the scalar case one cannot have f\/irst order dif\/ferential relations for the MOPRL, such as those given in~\eqref{ode1stOrder}. The general matrix case is much richer and complex, as  will be illustrated with some examples in the next section.

We f\/inish this section with a class of second-order dif\/ferential equations satisf\/ied by the MOPRL  $(\widehat{\bm P}_n)_n$. As a consequence of the Proposition~\ref{propladdX} we have the following
\begin{Proposition}\label{sodeq1X}
Under assumptions~\eqref{Ww} and~\eqref{ode:G}, with $\bm G(z)$ a polynomial, the MOPRL $(\widehat{\bm P}_n)_n$ satisfy the following second-order differential equation
\begin{gather}\label{seqex2}
    \widehat{\bm P}_n''+2\widehat{\bm P}_n' \bm G +\widehat{\bm P}_n\big(\bm G'+\bm G^2\big)
   +\mathcal{M}_n\widehat{\bm P}_n'+\mathcal{N}_n\widehat{\bm P}_n+\mathcal{M}_n\widehat{\bm P}_n \bm G=\bm 0,
\end{gather}
where $\mathcal{A}_n(z)=\mathcal{A}_n(z; \bm G)$,  $\mathcal{B}_n(z)=\mathcal{B}_n(z; \bm G)$,
\begin{gather*}%\label{MnzX}
   \mathcal{M}_n(z)= \mathcal{M}_n(z; \bm G)\\
   \phantom{\mathcal{M}_n(z)}{} =-(\mathcal{A}_n^{*}(z))'\mathcal{A}_n^{-*}(z)+\mathcal{B}_n(z)-\mathcal{A}_n^{*}(z)(z\bm I-\bm\alpha_n)+\mathcal{A}_n^{*}(z)\mathcal{B}_{n+1}(z)\mathcal{A}_n^{-*}(z),
\end{gather*}
and
\begin{gather*}\label{NnzX}
   \mathcal{N}_n(z)=\mathcal{N}_n(z; \bm G)=\mathcal{M}_n(z)\mathcal{B}_n(z)-\mathcal{B}_n^2(z)+\mathcal{B}_n'(z)+\mathcal{A}_n^{*}(z)\bm\beta_n\mathcal{A}_{n-1}^{*}(z),
\end{gather*}
provided that the inverse of $\mathcal{A}_n (z;\bm G)$ exists for $z\in\mathbb{C}$. Here, again, $-*$ denotes the conjugate transpose of the inverse.
\end{Proposition}

\begin{proof}
Dif\/ferentiate the lowering operator (\ref{lowop1X}) and substitute the raising operator (\ref{raiop1X}) evaluated at $n-1$.
\end{proof}

\begin{Remark}
Note that we can easily obtain another second-order dif\/ferential operator satisf\/ied by the orthogonal polynomials reversing the ladder operators. This new dif\/ferential equation needs not be in principle the same as~(\ref{seqex2}). Nevertheless, it is straightforward to see that both equations are equivalent using the compatibility conditions~(\ref{String11X}) and~(\ref{String21X}).
\end{Remark}

The equation \eqref{seqex2} does not have the form of the right hand side dif\/ferential operator considered for instance in \cite{MR2039133}, due to the terms $\mathcal{M}_n\widehat{\bm P}_n'$, $\mathcal{N}_n\widehat{\bm P}_n$ and $\mathcal{M}_n\widehat{\bm P}_n \bm G$. In some cases, under additional assumptions on the weight,  \eqref{seqex2}  can be reduced further, as we will see in the following section.

\section{Illustrative examples} \label{sec:examples}

In this section we study a number of examples of weights $\bm W$ which are smooth and non-vanishing on the whole real line; this assumption simplif\/ies the Riemann--Hilbert formulation because in this case one does not consider the local conditions (Y4) (see Section \ref{sec:RHformulation}). Additionally, without loss of generality, we take $ \bm W(0)=\bm I$.

For convenience, we consider the weights of the form
\begin{gather}\label{examplesW}
\bm W(x)=e^{-2q(x)} \bm U(x) \bm U^* (x), \qquad x\in \RR,
\end{gather}
where $q$ is a scalar real-valued function, so that with the notation \eqref{Ww} and \eqref{ode:G}, we may take
\begin{gather} \label{def:TandGpolynomials}
\bm T(x)=e^{-q(x)} \bm U(x) \qquad \text{and} \qquad \bm G(x)= - q'(x)\bm I + \bm U'(x) \bm U^{-1}(x).
\end{gather}

We are interested in the case when $\bm G$ is a matrix polynomial. Hence, keeping up with previous hypotheses, we will assume that $q$ is a (scalar) polynomial of even degree with real coef\/f\/icients and a positive leading coef\/f\/icient, and  that $\bm U'(x) \bm U^{-1}(x)$ is a matrix polynomial.

Since our main goal here is to generate a set of examples (some new, some already well known), we restrict the degree of $q$ to either $2$ (the Hermite case) or $4$ (the Freud case), and the degree of $\bm U'(x) \bm U^{-1}(x)$ to at most $1$. We start by considering the case when  $\bm U'(x) \bm U^{-1}(x)$ is a monomial, i.e.\ either $\bm U'(x) \bm U^{-1}(x)=\bm A$ or $\bm U'(x) \bm U^{-1}(x)=2\bm Bx$ for constant matrices $\bm A, \bm B\in\mathbb{C}^{N\times N}$. Taking into account the linearity~\eqref{linearity} of the coef\/f\/icients $\mathcal A_n$ and $\mathcal B_n$ in~\eqref{AnzX},~\eqref{BnzX},  we can obtain  the dif\/ferential equation \eqref{DeqX} for the general case of $\bm U'(x) \bm U^{-1}(x)=\bm A+2\bm Bx$. However, f\/inding the corresponding orthogonality weight is more involved:  when $\bm A$ and $\bm B$ do not commute, solving $\bm U'(x) \bm U^{-1}(x)=\bm A+2\bm Bx$ is not straightforward.  In Section~\ref{sec:4.3} we will discuss some examples, which yield explicit expressions for $\bm U(x)$, related to this case.

Recently an example has been found in \cite{BCD} where a weight matrix supported in the real line is explicitly given (but not of the type \eqref{examplesW}), when  $\bm G(x)$ is a matrix polynomial of degree $N$ with in general non-commuting coef\/f\/icients.

\subsection{The Hermite case}\label{sec:Herm}

For $q(x)= x^2/2$, let us consider two cases, f\/irst $\bm U'(x) \bm U^{-1}(x)=\bm A$ and then $\bm U'(x) \bm U^{-1}(x)=2\bm B x$. We end up this Section by discussing brief\/ly the case of $\bm U'(x) \bm U^{-1}(x)=\bm A+2\bm B x$.

\subsubsection[$U'(x) U^{-1}(x)=A$]{$\boldsymbol{U'(x) U^{-1}(x)=A}$}\label{sec:HermC}

The dif\/ferential equation $\bm U'(x) \bm U^{-1}(x)=\bm A$, so that $\bm U(x)= e^{\bm Ax}$, $\bm T(x)=e^{-x^2/2}e^{\bm Ax}$, and the weight matrix~\eqref{examplesW} is given by
\begin{gather} \label{FreudWeight}
   \bm W(x)=e^{-x^2}e^{\bm Ax}e^{\bm A^*x},\qquad\bm
A\in\mathbb{C}^{N\times N},\qquad x\in\mathbb{R}.
\end{gather}
The matrix $\bm G$ has the form
\begin{gather*} %\label{G:particularCase}
\bm G(x)=-  x  \bm I + \bm A,
\end{gather*}
and according to \eqref{Rm1x},
\begin{gather*}
 \mathcal{A}_{n}(x; \bm G) =2 \bm I, \qquad \mathcal{B}_n(x; \bm G)= - \bm G(x).
\end{gather*}
The compatibility conditions of Proposition \ref{propcc} yield
\begin{gather}\label{Comp1s}
2(\bm\beta_{n+1}-\bm\beta_n)=\bm I+\bm A\bm\alpha_n-\bm\alpha_n\bm A,
\end{gather}
and
\begin{gather}\label{Comp2s}
\bm \alpha_n=\frac{1}{2}(\bm A+\bm\gamma_{n}^{-1}\bm A^*\bm\gamma_{n}).
\end{gather}
The lowering and raising operators from Proposition~\ref{propladdX} are reduced now to
\begin{gather}\label{LowOps}
   \widehat{\bm P}_n'(x)+\widehat{\bm P}_n(x)\bm A-\bm A\widehat{\bm
P}_n(x)=2\bm\beta_n\widehat{\bm P}_{n-1}(x),
\end{gather}
and
\begin{gather*}%\label{RaiOps}
   -\widehat{\bm P}_n'(x)+2x\widehat{\bm P}_n(x)+\bm A\widehat{\bm
P}_n(x)-\widehat{\bm P}_n(x)\bm A-2\bm\alpha_n\widehat{\bm
P}_n(x)=2\widehat{\bm P}_{n+1}(x).
\end{gather*}
Summing up the telescopic relation in~\eqref{Comp1s} (or comparing the~$\mathcal{O}(x^{n-1})$ term in~(\ref{LowOps})) gives
\begin{gather}\label{betas}
   \bm\beta_n=\frac{1}{2}(n\bm I+\bm a_{n,n-1}\bm A-\bm A\bm a_{n,n-1}).
\end{gather}
With the notation of Proposition \ref{sodeq1X},
\begin{gather*}
\mathcal{M}_n(x) =    2 \left( \bm\alpha_n- \bm A\right),  \qquad
\mathcal{N}_n (x) =  -2 \left( \bm\alpha_n- \bm A\right)\bm G(x)-\bm G^2 (x) + \bm I + 4 \bm \beta_n ,
\end{gather*}
so that the dif\/ferential equation \eqref{seqex2} for the monic polynomials  $\widehat{\bm P}_n$, orthogonal with respect to the weight \eqref{FreudWeight}, boils down to
\begin{gather}
   \widehat{\bm P}_n''(x) +2\widehat{\bm P}_n'(x)(\bm A-x\bm
I)+\widehat{\bm P}_n(x)( \bm A^2 -2x\bm A)   \nonumber\\
\qquad{}= (-2x \bm A + \bm A^2 -4\bm\beta_n)\widehat{\bm P}_n(x) +2(\bm
A-\bm\alpha_n)(\widehat{\bm P}_n'(x)   +\widehat{\bm P}_n(x)\bm A-\bm
A\widehat{\bm P}_n(x)).\label{Secs}
\end{gather}

These formulas hold for any constant matrix $\bm A$, and we cannot expect important simplif\/ications without narrowing the class of the weights further. This can be done assuming in addition that the hypotheses of Proposition~\ref{thm:nonunique} hold. This, as it was shown in~\cite{MR2039133, MR2152234}, imposes additional constraints on the weight $\bm W$. Since the construction is described in detail in~\cite{MR2039133, MR2152234}, the exposition in this part will be rather sketchy.

The matrix $\bm H$ from \eqref{def:H}, restricted to $\RR$, can be written as
\begin{gather}
\bm H(x)=\bm T(x)\bm\chi \bm T^{-1}(x)=e^{\bm Ax}\bm\chi e^{-\bm
Ax}\nonumber\\
\phantom{\bm H(x)}{} =\bm\chi+\mbox{ad}_{\bm A}(\bm \chi)x+\mbox{ad}^2_{\bm
A}(\bm\chi)\frac{x^2}{2}+\cdots=\sum_{k\geq0}\mbox{ad}^k_{\bm
A}(\bm\chi)\frac{x^k}{k!}.\label{tht}
\end{gather}
Here $\bm\chi(x)=\bm S'(x) \bm S^* (x)$ is skew-Hermitian on $\RR$,  $\mbox{ad}_{\bm A}$ is the commutator given by $\mbox{ad}_{\bm A}(\bm\chi)=\bm A\bm\chi-\bm\chi \bm A$, and we def\/ine recursively
\begin{gather*}
\mbox{ad}^0_{\bm A}(\bm\chi)=\bm\chi, \qquad \mbox{ad}^{n+1}_{\bm A}(\bm\chi)=\mbox{ad}_{\bm
A}(\mbox{ad}^n_{\bm A}(\bm\chi))  \qquad \text{for}\quad  n\geq1.
\end{gather*}

The simplest situation obtains when the right hand side in (\ref{tht}) is constant; as it was shown in Lemma~2.4
of~\cite{MR2152234}, this assumption yields  $\bm \chi=ia\bm I$ for certain $a\in\mathbb{R}$ (which is the case discussed in Proposition~\ref{propsc}), when there are no new ladder
operators. Consequently, for the f\/irst non-trivial situation we must assume that~(\ref{tht}) is a matrix polynomial of degree at least
one. We consider two situations that yield degree exactly one (see~\cite{MR2039133, MR2152234} for motivations and further details).

Let us def\/ine a nilpotent matrix of the form
\begin{gather}\label{A}
     \bm L =\sum_{k=1}^{N-1}\nu_k\bm
E_{k,k+1},\qquad\nu_k\in\mathbb{C}\setminus\{0\},
\end{gather}
where $\bm E_{ij}$ is a matrix with 1 at entry $(i,j)$ and 0
elsewhere, and a diagonal matrix
\begin{gather}\label{JJ}
\bm J=\sum_{k=1}^N(N-k)\bm E_{k,k}.
\end{gather}

For the f\/irst non-trivial example we assume that $\bm A=\bm L$ and $\bm\chi=i\bm J$, so that
\begin{gather}\label{adCondition1}
\mbox{ad}_{\bm A}(\bm\chi)=-\bm A
\end{gather}
and  $\mbox{ad}^2_{\bm A}(\bm\chi)=\bm 0$.
Thus,
$\bm H(x)=e^{  {{\bm A}}x}\bm\chi e^{-  {{\bm A}}x}=i(\bm J-  {{\bm A}}x)$, and by   \eqref{Rm1x},
\begin{gather*}%\label{AnzXexH}
    \mathcal{A}_n(x; \bm H) =i(-  {{\bm A}}^*+\bm\gamma_n  {{\bm
A}}\bm\gamma_n^{-1})=2i(\bm\alpha_n^*-  {{\bm A}}^*), \\
  \mathcal{B}_n(x; \bm H) =i(  {{\bm A}}x-\bm J+\bm a_{n,n-1}  {{\bm
A}}-  {{\bm A}}\bm a_{n,n-1})=i(  {{\bm A}}x-\bm J+2\bm \beta_n-n\bm I).
\end{gather*}

From the  compatibility conditions of Proposition \ref{propcc} we get
\begin{gather}
\bm J\bm\alpha_n-\bm\alpha_n\bm J+\bm\alpha_n=  {{\bm
A}}+\frac{1}{2}\big(  {{\bm A}}^2\bm\alpha_n-\bm\alpha_n  {{\bm A}}^2\big) \qquad \text{and} \nonumber\\
 \bm J-\bm\gamma_{n}^{-1}\bm J\bm\gamma_{n}=  {{\bm
A}}\bm\alpha_n+\bm\alpha_n  {{\bm A}}-2\bm\alpha_n^2,\label{Comp2sH}
\end{gather}
where we have used relations (\ref{Comp1s}),
(\ref{Comp2s}), (\ref{betas}) and{\samepage
\begin{gather*}%\label{relbuena}
(\bm A-2\bm\alpha_n)\bm\beta_n=\bm\beta_n(\bm A-2\bm \alpha_{n-1}),
\end{gather*}
which can be easily proved using (\ref{betanttrr}) and (\ref{Comp2s}).}

Thus, from Corollary \ref{coro:additional}, (\ref{Comp1s}) and (\ref{Comp2sH}), the lowering and raising operators  (of the $0$-th order) are
\begin{gather*}%\label{LowOpsH}
\widehat{\bm P}_n(x)\bm J-\bm J\widehat{\bm P}_n(x)-x(\widehat{\bm
P}_n(x)  {{\bm A}}-  {{\bm A}}\widehat{\bm
P}_n(x))+2\bm\beta_n\widehat{\bm P}_n(x)-n\widehat{\bm P}_n(x)\\
\qquad{} =2(  {{\bm
A}}-\bm\alpha_n)\bm\beta_n\widehat{\bm P}_{n-1}(x),
\end{gather*}
and
\begin{gather*}%\label{RaiOpsH}
 \widehat{\bm P}_n(x)(\bm J-x  {{\bm A}})-\bm\gamma_{n}^{-1}(\bm
J-x  {{\bm A}}^*)\bm\gamma_{n}\widehat{\bm P}_n(x)+2\bm\beta_{n+1}
\widehat{\bm P}_n(x)-(n+1)\widehat{\bm P}_n(x)\\
\qquad{} =2(\bm \alpha_n-  {{\bm
A}})\widehat{\bm P}_{n+1}(x),
\end{gather*}
respectively, which yields by \eqref{ode1stOrder} the f\/irst order relation
\begin{gather}
   (  {{\bm A}}-\bm\alpha_n)\widehat{\bm P}_n'(x)+(  {{\bm
A}}-\bm\alpha_n+x\bm I)(\widehat{\bm P}_n(x)  {{\bm A}}-  {{\bm
A}}\widehat{\bm P}_n(x))-2\bm\beta_n\widehat{\bm P}_n(x)\nonumber\\
\qquad{} =\widehat{\bm
P}_n(x)\bm J-\bm J\widehat{\bm P}_n(x)-n\widehat{\bm P}_n(x).\label{fodeq}
\end{gather}
These identities hold in addition to the second order equation obtained in (\ref{Secs}), valid as we recall, for any $\bm A$. As far as we are aware of, these relations are new. We use them to simplify the dif\/ferential equation remarkably: multiplying  (\ref{fodeq}) by 2 and plugging it into   (\ref{Secs}) gives us for $\bm A= \bm L$,
\begin{gather*}%\label{SecOrd}
   \widehat{\bm P}_n''(x)+2 \widehat{\bm P}_n'(x)( {\bm A}-x\bm
I)+\widehat{\bm P}_n(x)\big( {\bm A}^2-2\bm J\big)=\big({-}2n\bm I+ {\bm
A}^2-2\bm J\big)\widehat{\bm P}_n(x),
\end{gather*}
which is a linear second-order dif\/ferential equation with coef\/f\/icients in the right hand side independent on $n$ (a.k.a.\ Sturm--Liouville equation with polynomial coef\/f\/icients), studied by Dur\'{a}n and Gr\"{u}nbaum in~\cite{MR2039133}.

In the second case we take
\begin{gather*}\label{BB}
    \bm A=  \bm L (\bm I+ \bm L)^{-1}=\sum_{j=1}^{N-1}(-1)^{j-1}  \bm L^j, \qquad \bm\chi=i\bm J,
\end{gather*}
where $ \bm L$ is given in (\ref{A}) and $\bm J$ in \eqref{JJ}. Consequently,
\begin{gather}\label{adCondition2}
\mbox{ad}_{\bm A}(\bm\chi)=-\bm A+\bm A^2, \qquad \mbox{ad}^2_{\bm A}(\bm\chi)=\bm 0,
\end{gather}
and
\begin{gather*}\bm H(x)=e^{ {{\bm A}}x}\bm\chi e^{- {{\bm A}}x}=i\big(\bm J-\big( {{\bm A}}- {{\bm A}}^2\big)x\big).
\end{gather*}
In this case,
\begin{gather*}%\label{AnzXexH2}
   \mathcal{A}_n(x; \bm H)=i\big({-} {{\bm A}}^*+( {{\bm A}}^*)^2+\bm\gamma_n\big( {{\bm A}}- {{\bm A}}^2\big)\bm\gamma_n^{-1}\big)\\
   \hphantom{\mathcal{A}_n(x; \bm H)}{} =2i(\bm
\alpha_n^*- {{\bm A}}^*-(\bm \alpha_n^*- {{\bm A}}^*)\bm \alpha_n^*-\bm
\alpha_n^*(\bm \alpha_n^*- {{\bm A}}^*))
\end{gather*}
and
\begin{gather*}%\label{BnzXexH2}
  \mathcal{B}_n(x; \bm H) =i\big(\big( {{\bm A}}- {{\bm A}}^2\big)x-\bm J+\bm a_{n,n-1}\big( {{\bm A}}- {{\bm A}}^2\big)-
  \big( {{\bm A}}- {{\bm A}}^2\big)\bm a_{n,n-1}\big) \\
 \phantom{\mathcal{B}_n(x; \bm H)}{}
  =i\big(\big( {{\bm A}}- {{\bm A}}^2\big)x-\bm J+2\bm \beta_n-n\bm I-2( {{\bm A}}\bm
\beta_n+\bm \beta_n {{\bm A}})+2n {{\bm A}}\big).
\end{gather*}
As to new results, we get the lowering operator,
\begin{gather*}
\widehat{\bm P}_n(x)\big(\bm J- x\big(  \bm A-  \bm A^2\big)\big)-\big(\bm J-x\big(  \bm A-  \bm A^2\big)\big)\widehat{\bm P}_n(x)\\
\qquad{} =(-2\bm \beta_n+n(\bm I-2  \bm A)+2(  \bm A\bm
\beta_n+\bm \beta_n  \bm A))\widehat{\bm P}_n(x)\\
\qquad\quad{}  -2(\bm \alpha_n-  \bm A-(\bm \alpha_n-  \bm A)\bm \alpha_n-\bm \alpha_n(\bm
\alpha_n-  \bm A))\bm \beta_n\widehat{\bm P}_{n-1}(x),
\end{gather*}
the raising operator
\begin{gather*}
\widehat{\bm P}_n(x)\big(\bm J- x\big(  \bm A-  \bm A^2\big)\big)-
\big(\bm J-x\big(  \bm A-  \bm A^2\big)\big)\widehat{\bm P}_n(x)\\
\qquad{} =(-2\bm \beta_n+n(\bm I-2  \bm A)+2(  \bm A\bm
\beta_n+\bm \beta_n  \bm A))\widehat{\bm P}_n(x)\\
\qquad\quad{} -2(\bm \alpha_n-  \bm A-(\bm \alpha_n-  \bm A)\bm \alpha_n-\bm \alpha_n(\bm
\alpha_n-  \bm A))((x\bm I-\bm \alpha_n)\widehat{\bm P}_{n}(x)-\widehat{\bm
P}_{n+1}(x)),
\end{gather*}
and the f\/irst-order dif\/ferential equation (\ref{ode1stOrder}):
\begin{gather*}%\label{fodeq2}
   (  \bm A-\bm \alpha_n- (  \bm A-\bm \alpha_n)\bm\alpha_n-\bm \alpha_n(  \bm A-\bm \alpha_n))\widehat{\bm P}_n'(x)\\
   \qquad{} =\widehat{\bm P}_n(x)\bm J-\bm
J\widehat{\bm P}_n(x)-x\big(\widehat{\bm P}_n(x)  \bm A^2-  \bm A^2\widehat{\bm P}_n(x)\big)\\
\qquad\quad{} +(x\bm I-  \bm A+\bm \alpha_n+(  \bm A-\bm
\alpha_n)\bm\alpha_n+\bm \alpha_n(  \bm A-\bm \alpha_n))(\widehat{\bm
P}_n(x)  \bm A-  \bm A\widehat{\bm P}_n(x))\\
 \qquad\quad{}+(2\bm \beta_n+n(2  \bm A-\bm I)-2(  \bm A\bm \beta_n+\bm
\beta_n  \bm A))\widehat{\bm P}_n(x).
\end{gather*}
Multiplying this equation by 2 and plugging it into the second-order
dif\/ferential equation (\ref{Secs}) gives
\begin{gather*}
   \widehat{\bm P}_n''(x) +2\widehat{\bm P}_n'(x)(  \bm A-x\bm
I)+\widehat{\bm P}_n(x)\big(  \bm A^2-2x  \bm A^2-2\bm J\big)=
\big(  \bm A^2-2x  \bm A^2-2\bm J\big)\widehat{\bm P}_n(x)\\
\qquad{}+(2n(2  \bm A-\bm I)-4(  \bm A\bm \beta_n+\bm \beta_n  \bm A)+2((\bm
\alpha_n-  \bm A)\bm \alpha_n+\bm \alpha_n(\bm \alpha_n-  \bm A))  \bm A)\widehat{\bm P}_n(x)\\
\qquad{} -2((\bm \alpha_n-  \bm A)\bm \alpha_n+\bm \alpha_n(\bm \alpha_n-  \bm A))(\widehat{\bm P}_n'(x)+\widehat{\bm P}_n(x)  \bm A).
\end{gather*}
This dif\/ferential equation is not of Sturm--Liouville type considered by Dur\'{a}n and Gr\"{u}nbaum in~\cite{MR2039133}, but nevertheless
we have been able to give a number of dif\/ferential equations of f\/irst and second order satisf\/ied by MOPRL with respect to a weight matrix that have not been considered up to this point.

\subsubsection[$U'(x) U^{-1}(x)=2 Bx$]{$\boldsymbol{U'(x) U^{-1}(x)=2 Bx}$}

The weight matrix \eqref{examplesW} is now given by
\begin{gather*}\label{WA2}
   \bm W(x)=e^{-x^2}e^{\bm Bx^2}e^{\bm B^*x^2},\qquad\bm
B\in\mathbb{C}^{N\times N},\qquad x\in\mathbb{R},
\end{gather*}
where we assume that $\bm B$ is chosen such that all the moments exist. By \eqref{def:TandGpolynomials},
\begin{gather*}
\bm T(x)=e^{-x^2/2}e^{\bm Bx^2} \qquad \text{and}\qquad \bm G(x)=(2\bm B-\bm I)x.
\end{gather*}
The weight matrix is an even matrix function, so that all moments of odd order vanish. As a~consequence,  $\bm a_{n,n-1}=\bm 0$ and $\bm \alpha_n=\bm 0$.
Then
\begin{gather*}
 \mathcal{A}_{n}(x; \bm G) =2\big(\bm I-\bm B^*-\bm \gamma_{n}\bm B\bm
\gamma_{n}^{-1}\big), \qquad \mathcal{B}_n(x; \bm G)= (\bm I-2\bm B)x,
\end{gather*}
and there will be only one compatibility condition:
\begin{gather*}%\label{Comp1s2}
2\big(\bm I-\bm B-\bm \gamma_{n+1}^{-1}\bm B^*\bm \gamma_{n+1}\big)\bm
\beta_{n+1}-2\bm \beta_n\big(\bm I-\bm B-\bm \gamma_{n-1}^{-1}\bm B^*\bm
\gamma_{n-1}\big)=\bm I.
\end{gather*}
From (\ref{LowOps2}) we get the explicit expression for $\bm \beta_n$ via
\begin{gather*}%\label{betas2}
   2\big(\bm I-\bm B-\bm \gamma_{n}^{-1}\bm B^*\bm \gamma_{n}\big)\bm
\beta_n=n\bm I+2(\bm a_{n,n-2}\bm B-\bm B\bm a_{n,n-2}).
\end{gather*}
The lowering and raising operators are
\begin{gather}\label{LowOps2}
   \widehat{\bm P}_n'(x)+2x(\widehat{\bm P}_n(x)\bm B-\bm B\widehat{\bm
P}_n(x))=2\big(\bm I-\bm B-\bm \gamma_{n}^{-1}\bm B^*\bm \gamma_{n}\big)\bm
\beta_n\widehat{\bm P}_{n-1}(x),
\end{gather}
and
\begin{gather*}%\label{RaiOps2}
   \widehat{\bm P}_n'(x)+2x\big(\widehat{\bm P}_n(x)\bm B-\bm B\widehat{\bm
P}_n(x)\big)=2\big(\bm I-\bm B-\bm \gamma_{n}^{-1}\bm B^*\bm
\gamma_{n}\big)\big(x\widehat{\bm P}_{n}(x)-\widehat{\bm P}_{n+1}(x)\big),
\end{gather*}
respectively. This yields  the
following second-order dif\/ferential equation
\begin{gather*}%\label{Secs2}
   \widehat{\bm P}_n''(x) +2x\widehat{\bm P}_n'(x)(2\bm B-\bm
I)+2x\big(\bm\gamma_n^{-1}\bm B^*\bm\gamma_n-\bm L_n\big)\widehat{\bm
P}_n'(x)+4x^2\big(\widehat{\bm P}_n(x)\bm B^2-\bm B^2\widehat{\bm
P}_n(x)\big)\\
 \qquad{}+4\bm K_n\widehat{\bm P}_n(x) +\big(\big(2-4x^2\big)\bm I+4x^2\big(\bm\gamma_n^{-1}\bm B^*\bm\gamma_n-\bm
L_n\big)\big)\big(\widehat{\bm P}_n(x)\bm B-\bm B\widehat{\bm P}_n(x)\big)=\bm 0,
\end{gather*}
where
\begin{gather*}%\label{Ln}
\bm L_n=\big(\bm I-\bm B-\bm \gamma_{n}^{-1}\bm B^*\bm \gamma_{n}\big)\bm B\big(\bm
I-\bm B-\bm \gamma_{n}^{-1}\bm B^*\bm \gamma_{n}\big)^{-1},
\end{gather*}
and
\begin{gather*}%\label{Kn}
\bm K_n=\big(\bm I-\bm B-\bm \gamma_{n}^{-1}\bm B^*\bm
\gamma_{n}\big)\bm\beta_n\big(\bm I-\bm B-\bm \gamma_{n-1}^{-1}\bm B^*\bm
\gamma_{n-1}\big).
\end{gather*}
These expressions are valid for any $\bm B$ (as long as all moments of the weight matrix exist). For further simplif\/ications we can assume again that the hypotheses of Proposition \ref{thm:nonunique} hold, and consider the cases given by the algebraic relations~\eqref{adCondition1} and~\eqref{adCondition2}. For instance, when~\eqref{adCondition2} holds, we obtain again a second-order dif\/ferential equation with coef\/f\/icients in the right hand side independent on~$n$, studied by Dur\'{a}n and Gr\"{u}nbaum in~\cite{MR2039133}. The computations are similar and will be omitted for the sake of brevity.

\subsubsection{Other cases}\label{sec:4.3}

Assume now that we have a linear combination of the previous two cases, i.e.
\begin{gather}
\label{exampleLinearComb}
\bm U'(x) \bm U^{-1}(x)=\bm A+2\bm Bx, \qquad \bm A, \bm B \in\mathbb{C}^{N\times N} , \qquad \text{with} \quad \bm U(0)=\bm I,
\end{gather}
in which case the matrix $\bm G(x)$ is given by
$   \bm G(x)= \bm A+(2\bm B-\bm I)x$.
This is all that is needed to calculate the coef\/f\/icients $\mathcal{A}_{n}(x; \bm G)$ and $\mathcal{B}_{n}(x; \bm G)$ and consequently the compatibility conditions, the ladder operators and the dif\/ferential relations.

However, for the corresponding orthogonality weight we need to solve \eqref{exampleLinearComb} explicitly, which may be non-trivial (unless $\bm A$ and $\bm B$ commute). In general, this  solution can be given in terms of the
\emph{time ordered} exponential
\begin{gather*}\label{WAOE}
   \bm U(x)= {:} \, e^{\int_0^x  (\bm A + \bm B s) ds}\, {:} , \qquad \bm A, \bm B \in\mathbb{C}^{N\times N},\qquad x\in\mathbb{R}.
\end{gather*}
This is obtained by rewriting the dif\/ferential equation as an integral one and ``solving'' it by iteration. In general, one cannot give an explicit expression for this inf\/inite sum.

An example of explicitly solvable non-trivial equation~\eqref{exampleLinearComb} can be found in \cite[Theorem~1.1]{Dtrio}, where non-commuting   $\bm A$ and $\bm B$ are constructed (as a linear combination of the matrices~$ {\bm L}$ and~$\bm J$ in~\eqref{A} and \eqref{JJ}, respectively) such that the corresponding $\bm U(x)=e^{( {\bm L}-v_0\bm J)x}e^{v_0\bm Jx}$, with $v_0$ any real number.  \cite{Dtrio}~contains also another example of a pair of non-commuting matrices  $\bm A$ and~$\bm B$ such that~$\bm U(x)$ is  a matrix polynomial of degree $N-1$, whose coef\/f\/icients can be obtained recursively in terms of~$\bm A$ and~$\bm B$ (see Theorem~A.1 therein).

An alternative way of generating examples is by considering weight matrices $\bm W$ of the form~\eqref{examplesW} with $\bm U(x)\in \mathbb P_m$ such that $\bm U(0)=\bm I$ and $\det \bm U(x)=1$ for $x\in \RR$. If we denote
\begin{gather*}
   \bm U(x)=\bm I+\sum_{k=1}^m\bm A_kx^k, \qquad  \bm U^{-1}(x)=\bm I+\sum_{k=1}^m\bm B_kx^k,
\end{gather*}
then the coef\/f\/icients $\bm A_k$ and $\bm B_k$ are connected by simple algebraic relations. For instance,
\begin{gather*}
  \bm B_1  =  -\bm A_1, \qquad
  \bm B_2  =  -\bm A_2+\bm A_1^2,\qquad
  \bm B_3  =  -\bm A_3+\bm A_2 \bm A_1+\bm A_1 \bm A_2-\bm A_1^3.
\end{gather*}
In consequence,
\begin{gather*}
 \bm U'(x) \bm U^{-1}(x)= \bm A_1+\big(2\bm A_2-\bm A_1^2\big)x+\big(3\bm A_3-2\bm A_2\bm A_1-\bm A_1\bm A_2+\bm A_1^3\big)x^2 \\
  \qquad{} + \big(4\bm A_4-3\bm A_3\bm A_1-\bm A_1\bm A_3+2\bm A_2\bm A_1^2+\bm A_1^2\bm A_2+\bm A_1\bm A_2\bm A_1-2\bm A_2^2-\bm A_1^4\big)x^3+\cdots.
\end{gather*}

Let us discuss some cases, depending on the degree $m$ of the matrix polynomial $\bm U(x)$.
In the f\/irst non-trivial example, when $m=1$, we have $\bm U(x)=\bm I+\bm A_1x$ and $\bm U^{-1}(x)=\bm I-\bm A_1x$ with $\bm A_1^2=\bm 0$.  Therefore, $\bm U'(x)\bm U^{-1}(x)=\bm A_1$, which was considered already in Section~\ref{sec:HermC}.

Let $m=2$, so that $\bm U(x)=\bm I+\bm A_1x+\bm A_2x^2$ and $\bm U^{-1}(x)=\bm I-\bm A_1x+(\bm A_1^2-\bm A_2)x^2$. The condition $\bm U(x)\bm U^{-1}(x)=\bm I$ yields the following algebraic relations between the coef\/f\/icients $\bm A_1$ and $\bm A_2$:
\begin{gather*}
    \bm A_2 \bm A_1+\bm A_1 \bm A_2=\bm A_1^3,\qquad \bm A_2\big({-}\bm A_2+\bm A_1^2\big)=\bm 0.
\end{gather*}
Algebraic manipulations of these two equations give additionally $\bm A_2 \bm A_1\bm A_2=\bm 0$ and $\bm A_2 \bm A_1^2=\bm A_1^2\bm A_2$. Therefore
\begin{gather*}
\bm G(x)=\bm A_1+\big(2\bm A_2-\bm A_1^2-\bm I\big)x-\bm A_2\bm A_1x^2.
\end{gather*}
The simplest case when $\bm G(x)$ is a matrix polynomial of degree one gives another algebraic equation, $\bm A_2\bm A_1=\bm 0$. This immediately implies that $\bm A_2^2=\bm 0$ and $\bm A_1^4=\bm 0$, i.e.\ once again resulting in very strong algebraic conditions on the coef\/f\/icients. In order to f\/ind an interesting example (when $\mbox{ad}_{\bm A_1}(\bm A_2)\neq\bm 0$) we have to go at least to the dimension $N=4$; for instance,
\begin{gather*}
\bm A_1=\begin{pmatrix}
  0 & a_{12} & a_{13} & a_{14} \\
  0 & 0 & a_{23} & a_{24} \\
  0 & 0 & 0 & a_{34} \\
  0 & 0 & 0 & 0 \\
\end{pmatrix},\qquad\bm A_2=\begin{pmatrix}
  0 & 0 & 0 & b_{14} \\
  0 & 0 & 0 & a_{23}a_{34} \\
  0 & 0 & 0 & 0 \\
  0 & 0 & 0 & 0 \\
\end{pmatrix},\\
\mbox{ad}_{\bm A_1}(\bm A_2)=\begin{pmatrix}
  0 & 0 & 0 & a_{12}a_{23}a_{34} \\
  0 & 0 & 0 & 0 \\
  0 & 0 & 0 & 0 \\
  0 & 0 & 0 & 0 \\
\end{pmatrix}.
\end{gather*}

Considering higher degree examples will give more and more algebraic relations among the coef\/f\/icients $\bm A_k$, which in general are dif\/f\/icult to solve.

In particular, in the examples that have appeared in the literature so far, $\bm U(x)$ is in general a matrix polynomial of degree depending on the size $N$.

\subsection{The Freud case}\label{sec:Freud}

Now we have $q(x)=\frac{x^4}{2}$, so the weight matrix is given by
\begin{gather*}
\bm W(x)=e^{-x^4} \bm U(x) \bm U^* (x), \qquad x\in \RR.
\end{gather*}
Again, we will consider two cases, f\/irst $\bm U'(x) \bm U^{-1}(x)=\bm A$ and then $\bm U'(x) \bm U^{-1}(x)=2\bm B x$.

\subsubsection[$U'(x) U^{-1}(x)= A$]{$\boldsymbol{U'(x) U^{-1}(x)= A}$}

The weight matrix \eqref{examplesW} in this case is given by
\begin{gather*}
   \bm W(x)=e^{-x^4}e^{\bm Ax}e^{\bm A^*x},\qquad\bm
A\in\mathbb{C}^{N\times N},\qquad x\in\mathbb{R}.
\end{gather*}
Then $\bm T(x)=e^{-x^4/2}e^{\bm Ax}$, $\bm G$ is given by
\begin{gather*} \label{FreudG}
\bm G(x)=  -2 x^3 \bm I  + \bm A.
\end{gather*}
Therefore using (\ref{AnzX}) and (\ref{BnzX}) we get
\begin{gather*}\label{AnzXexF1}
   \mathcal{A}_{n}(x; \bm G)=4\big(x^2\bm
I+\bm\alpha_nx+\bm\beta_{n+1}+\bm\beta_n+\bm\alpha_n^2\big)^*,
\end{gather*}
and
\begin{gather*}\label{BnzXexF1}
  \mathcal{B}_{n}(x; \bm G)=2x^3\bm
I+4(\bm\beta_nx+\bm\beta_n\bm\alpha_{n-1}+\bm\alpha_n\bm\beta_n)-\bm A.
\end{gather*}
These formulas are obtained by using the relation between the
coef\/f\/icients of the three-term recurrence equation~\eqref{ttrrm} and the
coef\/f\/icients of the MOPRL. In this case the coef\/f\/icients $\mathcal{A}_{n}(\cdot; \bm G)$ and $ \mathcal{B}_n(\cdot; \bm G)$
are matrix polynomials of degree~2 and~3, respectively.

From Proposition \ref{propcc} the compatibility conditions are
\begin{gather*}%\label{compatib1Freud}
\bm I+ \bm A\bm\alpha_n-\bm\alpha_n\bm
A-4(\bm\beta_{n+1}\bm\alpha_n+\bm\alpha_{n+1}\bm\beta_{n+1})\bm\alpha_n+4\bm\alpha_n(\bm\beta_{n}\bm\alpha_{n-1}+\bm\alpha_{n}\bm\beta_{n})\\
\qquad{} =4\big(\bm\beta_{n+2}+\bm\beta_{n+1}+\bm\alpha_{n+1}^2\big)\bm\beta_{n+1}-4\bm\beta_n(\bm\beta_{n}+\bm\beta_{n-1}+\bm\alpha_{n-1}^2),
\end{gather*}
and
\begin{gather*}
4\big(\big(\bm\beta_{n+1}+\bm\beta_n+\bm\alpha_n^2\big)\bm\alpha_n+\bm\alpha_n(\bm\beta_{n+1}+\bm\beta_n)+\bm\alpha_{n+1}\bm\beta_{n+1}+\bm\beta_n\bm\alpha_{n-1}
\big)=\bm
A+\bm\gamma_n^{-1}\bm A^*\bm\gamma_n.
\end{gather*}

The lowering operator is
\begin{gather*}
   \widehat{\bm P}_n'(x) +\widehat{\bm P}_n(x)\bm A-\bm A\widehat{\bm
P}_n(x)=-4(\bm\beta_nx+\bm\beta_n\bm\alpha_{n-1}+\bm\alpha_n\bm\beta_n)\widehat{\bm
P}_n(x)\\
 \qquad{} +4\big(x^2\bm
I+\bm\alpha_nx+\bm\beta_{n+1}+\bm\beta_n+\bm\alpha_n^2\big)\bm\beta_n\widehat{\bm
P}_{n-1}(x),
\end{gather*}
while the raising operator is
\begin{gather*}
   \widehat{\bm P}_n'(x) +\widehat{\bm P}_n(x)\bm A-\bm A\widehat{\bm
P}_n(x)\\
\qquad{} =4\big(x^3\bm
I+\bm\beta_{n+1}x-\big(\bm\beta_{n+1}+\bm\beta_n+\bm\alpha_n^2\big)\bm\alpha_n-\bm\beta_n\bm\alpha_{n-1}-\bm\alpha_n\bm\beta_n\big)\widehat{\bm
P}_n(x)\\
 \qquad\quad{} -4\big(x^2\bm
I+\bm\alpha_nx+\bm\beta_{n+1}+\bm\beta_n+\bm\alpha_n^2\big)\widehat{\bm
P}_{n+1}(x).
\end{gather*}

From the $\mathcal{O}(x^{n-1})$ term in the lowering operator we obtain
\begin{gather*}\label{betasF1}
n\bm I+\bm a_{n,n-1}\bm A-\bm A\bm
a_{n,n-1}
=4\big(\bm\beta_{n}\bm\beta_{n-1}+\bm\beta_{n}\bm\alpha_{n-1}^2+\bm\alpha_{n}\bm\beta_{n}\bm\alpha_{n-1}+\bm\beta_{n+1}\bm\beta_{n}
+\bm\beta_{n}^2+\bm\alpha_{n}^2\bm\beta_{n}\big).
\end{gather*}

An expression of the second-order dif\/ferential equation can be given using Proposition~\ref{thmXnp}, but in this case it is necessary to compute the inverse of $\mathcal{A}_{n}(x; \bm G)$, which is a matrix polynomial of degree~2 with a rational function determinant in general.

\subsubsection[$U'(x) U^{-1}(x)=2 Bx$]{$\boldsymbol{U'(x) U^{-1}(x)=2 Bx}$}

The weight matrix \eqref{examplesW} is given by
\begin{gather*}\label{WAF2}
   \bm W(x)=e^{-x^4}e^{\bm Bx^2}e^{\bm B^*x^2},\qquad\bm
B\in\mathbb{C}^{N\times N},\qquad x\in\mathbb{R},
\end{gather*}
we have $\bm T(x)=e^{-x^4/2}e^{\bm Bx^2}$ and
\begin{gather*}
\bm G(x)=-2x^3\bm I+2\bm Bx.
\end{gather*}
Again, this weight matrix is an even matrix function, and we conclude as above that $\bm a_{n,n-1}=\bm 0$
and $\bm \alpha_n=\bm 0$.   The coef\/f\/icients (\ref{AnzX}) and (\ref{BnzX}) are given by (using the
same approach as before)
\begin{gather*}
  \mathcal{A}_{n}(x; \bm G)=4\big(x^2\bm
I+\bm\beta_n^*+\bm\beta_{n+1}^*\big)-2\big(\bm B^*+\bm \gamma_{n}\bm B\bm
\gamma_{n}^{-1}\big),
\end{gather*}
and
\begin{gather*}
   \mathcal{B}_{n}(x; \bm G)=2\big(x^3\bm I+(2\bm\beta_n-\bm B)x\big).
\end{gather*}
Again, there is only one compatibility condition (\ref{String11X}):
\begin{gather*}
\bm
I=4((\bm\beta_{n+1}+\bm\beta_{n+2})\bm\beta_{n+1}-\bm\beta_n(\bm\beta_{n}+\bm\beta_{n-1}))\\
\hphantom{\bm I=}{}
-2\big(\bm
B\bm\beta_{n+1}-\bm\beta_{n}\bm B+\bm\gamma_{n+1}^{-1}\bm
B^*\bm\gamma_{n}-\bm\gamma_{n}^{-1}\bm B^*\bm\gamma_{n-1}\big)
\end{gather*}
with the lowering and raising operators
\begin{gather*}
   \widehat{\bm P}_n'(x) +2x\big(\widehat{\bm P}_n(x)\bm B-\bm B\widehat{\bm
P}_n(x)\big)+4x\bm\beta_n\widehat{\bm P}_n(x)\\
 \qquad{} =\big(4\big(x^2\bm I+\bm\beta_{n+1}+\bm\beta_n\big)-2\big(\bm B+\bm \gamma_{n}^{-1}\bm
B^*\bm \gamma_{n}\big)\big)\bm\beta_n\widehat{\bm P}_{n-1}(x),
\end{gather*}
and{\samepage
\begin{gather*}
   \widehat{\bm P}_n'(x) +2x\big(\widehat{\bm P}_n(x)\bm B-\bm B\widehat{\bm
P}_n(x)\big)=\big(4x^3\bm I+2\big(2\bm\beta_{n+1}-\bm B-\bm\gamma_n^{-1}\bm
B^*\bm\gamma_n\big)x\big)\widehat{\bm P}_n(x)\\
\qquad{}+ \big({-}4\big(x^2\bm I+\bm\beta_{n+1}+\bm\beta_n\big)+2\big(\bm B+\bm \gamma_{n}^{-1}\bm
B^*\bm \gamma_{n}\big)\big)\widehat{\bm P}_{n+1}(x),
\end{gather*}
respectively.}

Finally, considering the $\mathcal{O}(x^{n-2})$ term in the lowering operator we obtain
\begin{gather}\label{betasF2}
n\bm I+2(\bm a_{n,n-2}\bm B-\bm B\bm
a_{n,n-2})=4\big(\bm\beta_{n}\bm\beta_{n-1}+\bm\beta_{n}^2+\bm\beta_{n+1}\bm\beta_{n}\big)-2\big(\bm
B+\bm \gamma_{n}^{-1}\bm B^*\bm \gamma_{n}\big)\bm\beta_n.
\end{gather}
This equation may be regarded of as a matrix version of a
\emph{discrete Painlev\'{e} equation}. Observe that in the scalar
situation ($\bm B=\bm 0$), the equation \eqref{betasF2} reduces to the well-known string equations
\begin{gather*}
n=4\beta_n(\beta_{n+1}+\beta_n+\beta_{n-1}),
\end{gather*}
studied in the context of discrete Painlev\'{e}
equations and orthogonal polynomials of Freud type in \cite{MR2451211} (see
formula~(19) therein).

\section{Final remarks} \label{sec:finalremarks}

This paper extends to the matrix case the methodology of derivation of dif\/ferential relations
for MOPRL using a RH approach. We obtain explicit formulas for the ladder operators of very general weight matrices factorized in the form $\bm W=\bm T\bm T^*$. We have also shown that in the matrix case there is an extra freedom absent in the scalar situation, which allows us to obtain a family of ladder operators, some of them of the $0$-th order, which does not occur in the scalar case. Furthermore, by combining appropriately the family of ladder operators we  get a family of second order dif\/ferential equations satisf\/ied by the MOPRL, some of them  of the $0$-th or f\/irst order (under some invertibility assumptions). This yields some new results even in the particular cases studied before in the literature.

In order to keep the size of this paper reasonable, we have restricted our attention to the weights supported on the whole real line. When $\supp (\bm W)=[0,+\infty)$ or $\supp (\bm W)=[-1,1]$, we can follow the ideas exposed in Section \ref{sec:transf}, except that now we have to assume that either $z\bm G(z)$ (for the case of the semi-axis) or $(1-z^2)\bm G(z)$ (for the case of the f\/inite interval) are matrix polynomials. For the discussion of these situations in the scalar case, where the dif\/ferential equations for the Laguerre and Jacobi polynomials are obtained, the interested reader is referred to  \cite[Chapter~22]{Ismail05}.

\subsection*{Acknowledgements}

The authors are grateful for the excellent job of the referees, whose suggestions and remarks improved the f\/inal text.

The research of the f\/irst author was supported in part by the Applied Math.\ Sciences subprogram of the Of\/f\/ice of Energy Research, USDOE, under Contract DE-AC03-76SF00098.

The work of the second author is partially supported by the research project MTM2009-12740-C03-02 from the Ministry of Science and Innovation of Spain and the European Regional Development Fund (ERDF), by Junta de Andaluc\'{i}a grant  FQM-262, by K.U.~Leuven research grant OT/04/21, and by Subprograma de estancias de movilidad posdoctoral en el extranjero, MICINN, ref.~-2008-0207.

The third author is supported in part by Junta de Andaluc\'{\i}a grant FQM-229, and by the research project MTM2008-06689-C02-01 from the Ministry of Science and Innovation of Spain and the European Regional Development Fund (ERDF).

Both the second and the third authors gratefully acknowledge also the support of Junta de Andaluc\'{\i}a via the Excellence Research Grant P09-FQM-4643.

\pdfbookmark[1]{References}{ref}
\LastPageEnding

\end{document}